\documentclass[ijoc,sglanonrev]{classe}
\RequirePackage{tgtermes}
\RequirePackage{newtxtext}
\RequirePackage{newtxmath}
\RequirePackage{bm}
\RequirePackage{endnotes}

\OneAndAHalfSpacedXII 

\usepackage{tikz}

\usepackage{pgfplots}

\usepackage{caption}
\usepackage[labelfont=sf]{subcaption}
\usepackage{makecell}
\usepackage{amsmath, amssymb}
\pgfplotsset{compat=1.18}
\usepgfplotslibrary{groupplots}
\usepackage{multirow}
\usepackage{adjustbox}
\usepackage{comment}

\newcommand{\SBO}{{\texttt{S-BO}}}
\newcommand{\WBO}{{\texttt{W-BO}}}
\newcommand{\IBO}{{\texttt{I-BO}}}
\newcommand{\SWBO}{{\texttt{SWI-BO}}}
\newcommand{\SWDBO}{{\texttt{SWDI-BO}}}
\newcommand{\NIBO}{{\texttt{NI-BO}}}
\newcommand{\TIBO}{{\texttt{T-I-BO}}}
\newcommand{\Snew}{{\texttt{S}}}
\newcommand{\W}{{\texttt{W}}}
\newcommand{\I}{{\texttt{I}}}
\newcommand{\MCMC}{\texttt{MCMC}}
\newcommand{\DE}{{\texttt{DE}}}
\newcommand{\SAA}{{\texttt{SAA}}}
\newcommand{\MILP}{{\texttt{MILP}}}
\newcommand{\KKT}{{\texttt{KKT}}}

\usepackage[ruled]{algorithm2e}
\SetKw{Return}{Return:}
\SetKw{Returns}{returns}
\SetKw{Output}{Output:}
\SetKw{Input}{Input:}
\SetKw{Stepfive}{Step 5:}
\SetKw{Stepfour}{Step 4:}
\SetKw{Stepthree}{Step 3:}
\SetKw{Steptwo}{Step 2:}
\SetKw{Stepone}{Step 1:}
\SetKw{StepA}{a)}
\SetKw{StepB}{b)}
\SetKwFor{While}{While}{do}{end}
\SetKwFor{ForEach}{For each}{do}{end}
\SetKwFor{For}{For}{do}{end}
\SetKwIF{If}{ElseIf}{Else}{If}{then}{Else if}{Else}{end}

\usepackage{booktabs}

\usepackage{natbib}
 \bibpunct[, ]{(}{)}{,}{a}{}{,}%
 \def\bibfont{\small}%
\EquationsNumberedThrough    

\TheoremsNumberedThrough     
\ECRepeatTheorems  %

\MANUSCRIPTNO{IJOC-0001-2024.00}

\begin{document}

\RUNAUTHOR{Maria Bazotte, Margarida Carvalho, and Thibaut Vidal}

\RUNTITLE{Intermediate Bilevel Optimization: Modeling Endogenous Follower Tie-Breaking Behavior}

\TITLE{Intermediate Bilevel Optimization: Modeling Endogenous Follower Tie-Breaking Behavior}

\ARTICLEAUTHORS{

\AUTHOR{Maria Bazotte}
\AFF{MAGI \& CIRRELT, Polytechnique Montr\'{e}al, Montr\'{e}al, QC H3T 1J4, Canada, \EMAIL{maria-carolina.bazotte-corgozinho@polymtl.ca}}

\AUTHOR{Margarida Carvalho}
\AFF{DIRO \& CIRRELT, Universit\'{e} de Montr\'{e}al, Montr\'{e}al, QC H3T 1J4, Canada, \EMAIL{carvalho@iro.umontreal.ca}}

\AUTHOR{Thibaut Vidal}
\AFF{MAGI \& CIRRELT, Polytechnique Montr\'{e}al, Montr\'{e}al, QC H3T 1J4, Canada, \EMAIL{thibaut.vidal@polymtl.ca}}

} 

\ABSTRACT{%
In bilevel optimization, optimistic and pessimistic follower behaviors are the most commonly used forms to define how the follower ties-breaks among multiple optimal solutions. In this work, we go beyond these extreme tie-breaking behaviors and investigate the intermediate bilevel optimization program ({\IBO}), where the follower’s selected optimal response is a decision-dependent random event, with a probability measure influenced by the leader’s decision. We formally introduce a class of such endogenous measures, including the special case of strong-weak decision-dependent {\IBO}. We reformulate the {\IBO} as a Transformed {\IBO} ({\TIBO}) with exogenous uncertainty by defining inverse and Markov-chain transformations, which represent the follower’s response as a function of the leader’s decision and exogenous randomness. We handle the {\TIBO}'s uncertainty via sample-average approximation ({\SAA}), and we propose tailored approaches for its {\SAA} program according to the chosen transformation. Computationally, our methods solve reasonable-sized instances efficiently and outperform the deterministic equivalent when available. Furthermore, experiments stress the critical need to accurately model follower tie-breaking behavior, particularly depending on its alignment with the leader's objective, as misspecification leads to suboptimal leader decisions.
}%

\FUNDING{NSERC Grant No. 2024-04051, SCALE-AI research chair.}


\maketitle

\section{Introduction}\label{sec:introduction}

Bilevel optimization models hierarchical decision-making~\citep{bracken1973mathematical}, requiring a defined tie-breaking rule when the follower's problem admits multiple optimal responses. While literature traditionally defaults to extreme optimistic (fully cooperative) or pessimistic (fully adversarial) follower behaviors~\citep[see, e.g.,][]{beck2023survey}, these assumptions often expose the leader to excessive risk or missed gains. In practice, followers frequently exhibit intermediate, nuanced levels of cooperation, such as in competitive settings~\citep{corpus2025study}, making it critical to utilize data or beliefs to model more generalized tie-breaking behaviors rather than relying on inappropriate extremes~\citep{beck2023survey,salas2023existence}.

To address the limitations of extreme assumptions, the intermediate approach models the leader's belief regarding the follower's tie-breaking behavior as a decision-dependent probability measure over their optimal response set~\citep{lina1996hierarchical}. Because this behavior may alter the leader's optimal solutions~\citep{lina1996hierarchical,jia2011new,zeng2020practical}, it requires careful modeling. Following~\citet{jia2016new}, we assume the follower behaves rationally: its cooperativeness monotonically decreases as the leader's decision yields worse optimal outcomes for the follower. This aligns the model with observed incentives and avoids pathological behaviors, such as followers rewarding unfavorable leader decisions.  

However, solving these intermediate bilevel models is methodologically challenging. Optimistic bilevel programs are already strongly NP-hard~\citep{hansen1992new}. The intermediate approach is at least as difficult, as even the particular case of a fixed convex combination of optimistic and pessimistic behaviors remains strongly NP-hard when optimistic and pessimistic solutions are known~\citep{lagos2023complexity}. Furthermore, the follower's response introduces Type 1 endogenous uncertainty~\citep{goel2006class}, where the leader's decision dictates both the probability measure and the support. This renders classical methods, like deterministic equivalent ({\DE}) formulations, highly complex and non-convex~\citep{dupacova2006optimization}. Since existing literature for the intermediate approach primarily focuses on existence results rather than formal modeling or tractable algorithms~\citep[see][]{beck2023survey,lina1996hierarchical,mallozzi2005oligopolistic,salas2023existence,cotrina2025lipschitzianity}, there is a critical need for formal endogenous measures and practical solution methods.

Nevertheless, the literature on the intermediate approach remains limited~\citep[see][]{beck2023survey}. To the best of our knowledge, only a few studies~\citep{lina1996hierarchical,mallozzi2005oligopolistic,salas2023existence,cotrina2025lipschitzianity} investigated the general intermediate bilevel framework, and these primarily establish existence results rather than formally defining decision-dependent probability measures for follower behavior or proposing practical solution methods. This highlights the need for formal decision-dependent measures that capture realistic follower beliefs, along with solution methods that address both the NP-hardness of bilevel problems and the added complexity of decision-dependent uncertain follower responses, which are the focus of this work.

\paragraph{Contributions.}

We introduce endogenous probability measures that capture different follower behaviors, along with a tractable methodology to solve the resulting {\IBO}. Our contributions are:
\begin{enumerate}
    \item In Section~\ref{sec:modeling}, we formalize endogenous probability measures for two cases: (a) dichotomous measures with support restricted to optimistic and pessimistic solutions, and (b) generalized measures with support possibly covering the entire follower's optimal set of responses.
    \item In Section~\ref{sec:methodology}, we propose a practical methodology to tackle the {\IBO}. We are the first to reformulate this complex program with endogenous uncertainty into a Transformed {\IBO} ({\TIBO}) with purely exogenous uncertainty by applying the paradigm of~\citet{bazotte2025solving}. We then build on the Sample Average Approximation ({\SAA}) framework to address the resulting uncertainty in the {\TIBO}. For dichotomous measures, we use an inverse transformation and adapt~\citet{zeng2020practical} to obtain an exact, correction-free two-level {\SAA} program. For generalized measures, we use a Markov chain transformation modeled implicitly via integer {\sf L}-shaped cuts and an appropriate Markov Chain Monte Carlo {\MCMC} simulation, for problems with binary leader decisions and convex follower problems.
    \item In Section~\ref{sec:experiments}, we conduct experiments on integer-linear upper-level problems without coupling constraints and a linear follower problem, adapted from the BOBILib library~\citep{BOBILib:2026}.
    \begin{enumerate}
        \item We show that, for dichotomous measures where the {\DE} program baseline approach is available, the proposed {\SAA} method is 1.8 times faster on average, while maintaining improved solution quality. Furthermore, solving the {\SAA} program of the {\IBO} with the proposed endogenous measures is considerably more challenging than solving the optimistic and pessimistic approaches, highlighting the impact of follower behavior with decision-dependent probabilities. 
        \item We demonstrate that accurately modeling follower behavior is critical. In particular, when the follower’s objective is more aligned to the leader’s, the follower can improve its own outcome through strategic behavior. Moreover, misspecifying follower behavior yields worse leader solutions in terms of the leader’s objective than those obtained under the actual behavior, underscoring the critical need for accurate modeling and the robustness of balanced tie-breaking behaviors.
    \end{enumerate}
\end{enumerate}

\section{Problem Definition and Related Works}\label{sec:problem-literature}

We focus on \emph{intermediate bilevel programs} ({\IBO}s). Let ${x \in \mathcal{X} \subseteq \mathbb{R}^{n_x}_{+}}$ denote the leader's decisions, and, for any given $x$, let ${y \in \mathcal{Y}(x) \subseteq \mathbb{R}^{n_y}_{+}}$ denote the follower's feasible decisions. Following prior work~\citep[see][]{beck2023survey}, we assume that ${\mathcal{X} = \{x \in X : G(x) \ge 0\}}$ is a nonempty compact set without coupling constraints, and that ${\mathcal{Y}(x) = \{y \in Y : g(x,y) \ge 0\}}$ is nonempty and compact for all $x$. The sets $X$ and $Y$ impose nonnegativity, integer, or binary restrictions. For a given $x \in \mathcal{X}$, the follower's optimal value function and $\varepsilon$-optimal response set are defined as:
\begin{equation*}
    \varphi(x) = \min_{y} \{ f(x,y) : y \in \mathcal{Y}(x) \}, \qquad
    \mathcal{S}(x,\varepsilon) = \{y \in \mathcal{Y}(x) : f(x,y) \leq \varphi(x) + \varepsilon \}, \quad \varepsilon \geq 0.
\end{equation*}
The $\varepsilon > 0$ enables assessing whether a follower may deviate from exact optimum to influence the leader’s decision toward outcomes that are ultimately more favorable to them. The {\IBO} is:
\begin{equation}\label{Prog:intermediate}
    \mbox{{\IBO}}: \Theta_{\I}({\mu}) = \min_{x\in \mathcal{X}} \left\{ H(x) + \mathbb{E}_{\boldsymbol{\hat{y}}\! \sim \! \mu_x}\! \left[ F(x,\boldsymbol{\hat{y}}) \right] \right\},
\end{equation}
where, for each leader decision $x \in \mathcal{X}$, $\mu_x$ is a probability measure supported on $\operatorname{supp}(\mu_x) \subseteq \mathcal{S}(x,\varepsilon)$ capturing the leader’s belief about the follower’s response behavior. We denote $\mu = \cup_{x \in \mathcal{X}} \{\mu_x\}$. Finally, we assume that the leader’s and follower’s objective and constraint functions ($H$, $F$, $G$, $f$, $g$) are well-defined and continuous on their respective domains.

We denote as $\hat{y} \in \mathcal{S}(x,\varepsilon)$ a realization of the endogenous random vector $\boldsymbol{\hat{y}}\! \sim \!\mu_x$. The expected follower objective is $\mathbb{E}_{\boldsymbol{\hat y}\!\sim\!\mu_x}[f(x,\boldsymbol{\hat y})]$, which may differ from $\varphi(x)$ when $\varepsilon >0$. Appendix~\ref{app:intermediate-discrete-continuous} details the {\IBO} for finite $\mathcal{S}(x,\varepsilon)$ or absolutely continuous $\mu_x$. Moreover, the {\IBO} is a stochastic program with decision-dependent uncertainty. The expectation $\mathbb{E}_{\boldsymbol{\hat{y}}\sim\mu_x}[F(x,\boldsymbol{\hat{y}})]$ generally lacks a closed-form expression and may be nonconvex in $x$, making {\DE} formulations and sampling highly complex, as both the distribution and its support $\mathcal{S}(x,\varepsilon)$ depend on $x$~\citep{dupacova2006optimization}.

\subsection{Related Works}\label{subsec:literature-background}

Table~\ref{tab:bilevel_measures} summarizes the two {\IBO} classes reviewed next: \emph{support-endogenous}, where measures are fixed probabilities with decision-dependent support, and \emph{fully endogenous}, where measures have both support and probabilities that are decision-dependent.

\paragraph{\textbf{Support-endogenous.}} The {\IBO} generalizes the optimistic, or strong, bilevel program ({\SBO}) and the pessimistic, or weak, bilevel program ({\WBO}):
\begin{equation*}
    \mbox{{\SBO}: } \Theta_{\I}(\mu_{\Snew}) = \min_{x \in \mathcal{X}} \left\{  H(x) +  \min_{\hat{y} \in \mathcal{S}(x,\varepsilon)} F(x,\hat{y}) \right\}, \quad \mbox{{\WBO}: } \Theta_{\I}(\mu_{\W}) = \min_{x \in \mathcal{X}} \left\{ H(x) + \max_{\hat{y} \in \mathcal{S}(x,\varepsilon)} F(x,\hat{y}) \right\},
\end{equation*}
whose leader values bound the {\IBO}~\citep{lina1996hierarchical}.

The {\IBO} was first investigated by~\citet{aboussoror1995strong} and~\citet{lina1996hierarchical}, the latter introducing the general {\IBO} in single-leader, single- and multi-follower settings, with both optimal and near-optimal follower response sets. They established existence results for the {\IBO} and illustrated that the optimistic and pessimistic solutions may differ substantially from intermediate ones. Subsequent works further studied the existence and structural properties of the {\IBO} in oligopolistic markets and multi-leader problems~\citep{mallozzi2005oligopolistic,salas2023existence}. More recently,~\citet{cotrina2025lipschitzianity} established sufficient conditions for Lipschitz continuity of the expected value function in stochastic programs with decision-dependent uncertainty and moving support, a class that includes the {\IBO}.

A commonly studied special case is the strong-weak {\IBO} ({\SWBO}), which uses a fixed probability $\beta \in [0,1]$ to select between the optimistic and pessimistic solutions:
\begin{equation*}
    \mbox{\SWBO}: \Theta_{\I} (\mu_\beta) = \min_{x\in \mathcal{X}} \left\{ H(x) +  \beta \! \min_{\hat{y} \in \mathcal{S}(x,\varepsilon)} F(x,\hat{y}) + \left[ 1 - \beta \right] \! \max_{\hat{y}\in \mathcal{S}(x,\varepsilon)} F(x,\hat{y}) \right\}.
\end{equation*}

The {\SWBO} was first studied by~\citet{aboussoror1995strong}, who established existence results in finite-dimensional spaces, later extended to infinite-dimensional settings by~\citet{aboussoror2017strong}.~\citet{cao2002partial,zheng2015new} and~\citet{zeng2020practical} studied the linear {\SWBO}, in which both upper- and lower-level problems are linear programs without upper-level coupling constraints.~\citet{cao2002partial} proposed penalty-based reformulations that reduce the problem to a classical linear bilevel program solvable via standard techniques such as {\KKT} reformulation, and illustrated numerically how varying $\beta$ yields optimistic, pessimistic, or intermediate leader decisions.~\citet{zheng2015new} introduced an exact penalty method and analyzed the leader’s optimal value as a function of $\beta$, deriving a procedure to identify its critical points without exhaustive enumeration.~\citet{zeng2020practical} proposed a relaxation-correction scheme that reformulates the trilevel structure of the {\SWBO} into a classical bilevel program with two followers, corresponding to the optimistic and pessimistic cases, solvable via standard methods. In their approach, the pessimistic follower’s optimality conditions are first relaxed and subsequently enforced with a correction step.

\paragraph{\textbf{Fully endogenous.}} Extending the {\SWBO}, the strong--weak decision-dependent {\IBO} ({\SWDBO}) assumes the cooperation level $\beta(x):\mathcal{X} \!\to\! [0,1]$ depends on the leader's decision:
\begin{equation}\label{Prog:swd}
    \mbox{\SWDBO}: \Theta_{\I}(\mu_{\beta(\cdot)}) = \min_{x\in \mathcal{X}} \left\{ H(x) + \beta \!\left( x \right)  \! \min_{\hat{y} \in \mathcal{S}(x,\varepsilon) } F(x,\hat{y}) + \left[ 1 - \beta \! \left( x \right) \right] \! \max_{\hat{y}\in \mathcal{S}(x,\varepsilon)} F(x,\hat{y}) \right\}.
\end{equation}

\citet{jia2011new} introduced the {\SWDBO}, established existence results under lower semicontinuity of the follower’s optimal response set, and provided examples where the {\SBO} and {\WBO} yield distinct optimal leader solutions from the {\SWDBO}.~\citet{jia2016new} modeled the cooperation level $\beta(x)$ as the follower’s satisfaction with the leader’s decision, assuming that greater satisfaction induces stronger cooperation. They proved existence results under mild assumptions and showed that, under certain conditions, the {\SWDBO} yields solutions with higher joint leader--follower satisfaction than the fixed-cooperation {\SWBO}. To our knowledge, these remain the only studies on the {\SWDBO}; research on this topic is still scarce~\citep[see][]{beck2023survey}, particularly regarding the practical modeling of the decision-dependent function $\beta(x)$, including its definition and structural properties, as well as solution methods.

The neutral-measure {\IBO} ({\NIBO}), with $\mu_x^{\text{Ntr}}$ the uniform distribution over $\mathcal{S}(x,\varepsilon)$ for each $x \in \mathcal{X}$ and $\mu_{\text{Ntr}} = \cup_{x \in \mathcal{X}} \{\mu_x^{\text{Ntr}} \}$, has been investigated by~\citet{salas2023existence}. The {\NIBO} is:
\begin{equation}
    \mbox{NI-BO: } \min_{x \in \mathcal{X}} \left\{ H(x) \!+\! \frac{1}{\lambda(\mathcal{S}(x,\varepsilon))} \int_{\hat{y} \in \mathcal{S}(x,\varepsilon)} \!\! F(x,\hat{y}) \, d \lambda(\hat{y}) \right\}, \; \mu_x^{\text{Ntr}}(A) \!=\! \frac{\lambda\left(A \cap \mathcal{S}(x,\varepsilon)\right)}{\lambda\left(\mathcal{S}(x,\varepsilon)\right)},  \forall\ A \subseteq \mathbb{R}^{n_y},  \nonumber
\end{equation}
where $\lambda(\cdot)$ is the Lebesgue measure. To our knowledge, the only solution procedure proposed for this setting is the differential evolution metaheuristic with inner Monte Carlo sampling by~\citet{salas2023existence} for continuous (non) linear upper-level problems with linear follower problems.
\begin{table}[t]
    \TABLE
    {Summary of {\IBO} classes: Support- and Full-endogenous classifications. \label{tab:bilevel_measures}}
    {
    \centering
    \begin{adjustbox}{width=0.68\textwidth,center}
    \centering
    \begin{tabular}{r|ccccc|c}
        \hline
        & \multicolumn{5}{c|}{$\longleftarrow$ \text{Less General \quad \quad \quad More General} $\longrightarrow$}  & Uniform \\
        & {\SBO} & {\WBO} & {\SWBO} & {\SWDBO}  & {\IBO} & {\NIBO}\\ \hline
        \textrm{Notation} $\mu$ & $\mu_{\Snew}$ & $\mu_{\W}$ & $\mu_\beta$ & $\mu_{\beta(\cdot)}$  & $\mu$ & ${\mu_{\text{Ntr}}}$\\
        $\operatorname{supp}(\mu_x)$ & $\lbrace \hat{y}^{\Snew}(x,\varepsilon) \rbrace$ & $\lbrace \hat{y}^{\W}(x,\varepsilon) \rbrace $ & $\{\hat{y}^{\Snew}(x,\varepsilon),\hat{y}^{\W}(x,\varepsilon)\}$ & $\{\hat{y}^{\Snew}(x,\varepsilon),\hat{y}^{\W}(x,\varepsilon)\}$ & $\mathcal{S}(x,\varepsilon)$ & $\mathcal{S}(x,\varepsilon)$ \\[-1ex]
        & {\tiny (Singleton)} & {\tiny (Singleton)} & {\tiny (Dichotomous)} & {\tiny (Dichotomous)} & {\tiny (Generalized)} & {\tiny (Generalized)} \\
        \textrm{Type} & \text{Support} & \text{Support} & \text{Support} & \text{Full} & \text{Full} & \text{Full} \\\hline 
   \end{tabular}
   \end{adjustbox}
   }
   {We define ${\hat{y}^{\Snew}(x,\varepsilon) \in \arg \! \min_{\hat{y}}\{F(x,\hat{y}): \hat{y} \in \mathcal{S}(x,\varepsilon)\}}$ and ${\hat{y}^{\W}(x,\varepsilon) \in \arg\!\max_{\hat{y}}\{F(x,\hat{y}): \hat{y} \in \mathcal{S}(x,\varepsilon)\}}$ as the optimistic and pessimistic follower $\varepsilon$-optimal solutions.}
\end{table}

\paragraph{\textbf{Endogenous uncertainty.}} While related literature explores multi-objective lower levels, side-payments, and adversarial deviations~\citep{jia2013new,zare2018class,alves2021new}, the core challenge of the {\IBO} stems from its Type 1 endogenous uncertainty, a setting that has received limited attention in prior work~\citep{li2021review}. Initial works have focused on solving the typically nonconvex {\DE} using standard solvers~\citep{hellemo2018decision}, extensions of stochastic gradient~\citep{homem2022simulation} and {\sf L}-shaped methods~\citep{pantuso2025shaped} for endogenous uncertainty, greedy algorithms~\citep{karaesmen2004overbooking}, and reformulations based on random variable transformations~\citep{bazotte2025solving}. In the context of the {\IBO}, fully endogenous {\IBO} measures remain largely unexplored, as they combine the challenges of endogenous uncertainty and bilevel optimization.

\section{Modeling Intermediate Decision-Dependent Follower Behavior}\label{sec:modeling}

Following~\citet{jia2016new}, we define full-endogenous measures where the leader's decision influences the follower's behavior. We assume that the follower becomes more cooperative as the leader’s decision yields a better outcome for them and more adversarial as it worsens. Let ${\varphi_{lb} := \min_{x' \in \mathcal{X}} \varphi(x')}$ and ${\varphi_{ub} := \max_{x' \in \mathcal{X}} \varphi(x')}$ be finite constants bounding  $\varphi(x)$ for all $x \in \mathcal{X}$. Moreover, for any realization $\hat{y}$ of $\boldsymbol{\hat{y}} \sim \mu_x$, the leader's reaction-dependent objective component satisfies $F_{\Snew}(x,\varepsilon) \leq F(x,\hat{y}) \leq F_{\W}(x,\varepsilon)$, where $F_{\Snew}(x,\varepsilon) = F\left(x, \hat{y}^{\Snew}(x,\varepsilon)\right)$ and $F_{\W}(x,\varepsilon) = F\left(x,\hat{y}^{\W}(x,\varepsilon)\right)$ are the leader optimistic and pessimistic values. Under the proposed measures, the follower becomes more cooperative when $\varphi(x)$ is close to $\varphi_{lb}$, assigning higher probability to leader-favorable $\varepsilon$-optimal responses with lower values of $F(x,\hat y)$, close to $F_{\Snew}(x,\varepsilon)$. As $\varphi(x)$ approaches $\varphi_{ub}$, the follower becomes less cooperative, assigning higher probability to responses with larger values of $F(x,\hat y)$, close to $F_{\W}(x,\varepsilon)$. We distinguish two classes of measures: dichotomous (Section~\ref{subsec:modeling-swd}), which assign positive probability only to the optimistic and pessimistic solutions, and generalized (Section~\ref{subsec:modeling-ibl}), which may assign positive probability to any near-optimal follower solution.

\subsection{Strong Weak Decision-Dependent Measures}\label{subsec:modeling-swd}

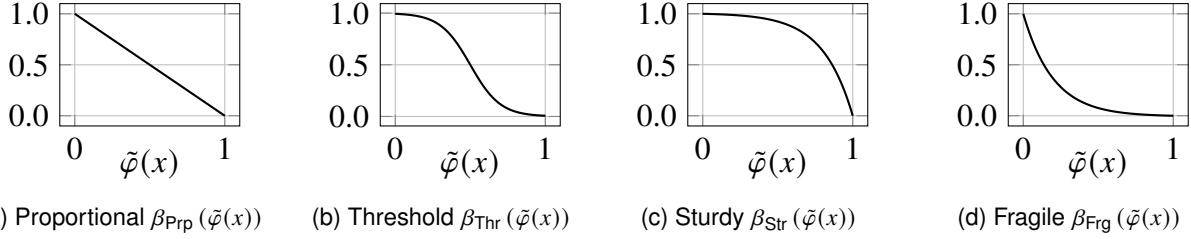
\begin{figure}[t]
    \caption{Proposed cooperation $\beta \left( \tilde{\varphi}(x) \right)$ functions over $\tilde{\varphi}(x) \in [0, 1]$.}
    \label{fig:beta-functions}
    \centering
    \vspace{0.1cm}
    \begin{subfigure}{0.24\textwidth} 
        \centering
        \begin{tikzpicture}
        \begin{axis}[
            width=\textwidth, height=3.2cm,
            xlabel={$\tilde{\varphi}(x)$},
            xtick={0,1}, xticklabels={0,1},
            ytick={0,0.5,1}, yticklabels={0.0,0.5,1.0},
            xlabel style={yshift=12pt},
            grid=both
        ]
        \addplot[domain=0:1, samples=100, thick, black]{1 - x}; 
        \end{axis}
        \end{tikzpicture}
        \caption{Proportional $\beta_{\text{Prp}} \left( \tilde{\varphi}(x) \right)$}
    \end{subfigure}
    \hfill
    \begin{subfigure}{0.24\textwidth}
        \centering
        \begin{tikzpicture}
        \begin{axis}[
            width=\textwidth, height=3.2cm,
            xlabel={$\tilde{\varphi}(x)$},
            xtick={0,1}, xticklabels={0,1},
            ytick={0,0.5,1}, yticklabels={0.0,0.5,1.0},
            xlabel style={yshift=12pt},
            grid=both
        ]
        \addplot[domain=0:1, samples=100, thick, black]{1/(1 + exp(10*(x - 0.5)))}; 
        \end{axis}
        \end{tikzpicture}
        \caption{Threshold $\beta_{\text{Thr}} \left( \tilde{\varphi}(x) \right)$}
    \end{subfigure}
    \begin{subfigure}{0.24\textwidth}
        \centering
        \begin{tikzpicture}
        \begin{axis}[
            width=\textwidth, height=3.2cm,
            xlabel={$\tilde{\varphi}(x)$},
            xtick={0,1}, xticklabels={0,1},
            ytick={0,0.5,1}, yticklabels={0.0,0.5,1.0},
            xlabel style={yshift=12pt},
            grid=both
        ]
        \addplot[domain=0:1, samples=100, thick, black]{(1 - exp(-5*(1-x)))/(1 - exp(-5))}; 
        \end{axis}
        \end{tikzpicture}
        \caption{Sturdy $\beta_{\text{Str}} \left( \tilde{\varphi}(x) \right)$}
    \end{subfigure}
    \hfill
    \begin{subfigure}{0.24\textwidth}
        \centering
        \begin{tikzpicture}
        \begin{axis}[
            width=\textwidth, height=3.2cm,
            xlabel={$\tilde{\varphi}(x)$},
            xtick={0,1}, xticklabels={0,1},
            ytick={0,0.5,1}, yticklabels={0.0,0.5,1.0},
            xlabel style={yshift=12pt},
            grid=both
        ]
        \addplot[domain=0:1, samples=100, thick, black]{1 - (1 - exp(-5*(x)))/(1 - exp(-5))}; 
        \end{axis}
        \end{tikzpicture}
        \caption{Fragile $\beta_{\text{Frg}} \left( \tilde{\varphi}(x) \right)$}
    \end{subfigure}
\end{figure}

In the {\SWDBO} model, we assume that $\beta(x)$ depends on $x$ solely through the follower’s optimal value $\varphi(x)$. To formalize this, we introduce the normalized optimal follower value $\tilde{\varphi}(x) \in [0,1]$ as:
\begin{equation*}
    \tilde{\varphi}(x) = \frac{\varphi(x) - \varphi_{lb}}{\varphi_{ub} - \varphi_{lb}} \in [0,1], \qquad \forall\ x \in \mathcal{X}, \qquad \text{ with } \, \varphi_{lb} < \varphi_{ub},
\end{equation*}
where $\tilde{\varphi}(x) = 0$ corresponds to the best possible outcome for the follower. We assume that:
\begin{assumption}\label{assump:mono}
    The cooperation function $\beta \left( \tilde{\varphi}(x) \right) \mapsto [0,1]$ is strictly decreasing in $\tilde{\varphi}(x)$.
\end{assumption}
We impose strict monotonicity to capture that every change in the follower’s optimal objective value $\varphi(x)$ alters the cooperation level, leaving no flat segments in $\beta \left( \tilde{\varphi}(x) \right)$. We consider four distinct behavioral forms: Proportional (gradual decrease), Threshold (sharp drop at a central threshold), Sturdy (highly cooperative until approaching $\varphi_{ub}$), and Fragile (rapid decline in cooperation once exceeding $\varphi_{lb}$). The mathematical definitions for these dichotomous measures--$\beta_{\text{Prp}}(\cdot)$, $\beta_{\text{Thr}}(\cdot)$, $\beta_{\text{Str}}(\cdot)$, $\beta_{\text{Str-p}}(\cdot)$, $\beta_{\text{Frg}}(\cdot)$, $\beta_{\text{Frg-p}}(\cdot)$--are summarized in Table~\ref{tab:cooperation_functions} (introduced in the next section alongside generalized measures). A illustration of these behaviors is provided in Figure~\ref{fig:beta-functions} and their impacts on optimal decisions are shown in Example~\ref{ex:behavior-swd} (Appendix~\ref{app:examples}). Finally, we denote ${\mathcal{B} = \{\beta_{\text{Prp}}(\cdot), \beta_{\text{Thr}}(\cdot), \beta_{\text{Str}}(\cdot),\beta_{\text{Str-p}}(\cdot),\beta_{\text{Frg}}(\cdot),\beta_{\text{Frg-p}}(\cdot)\}}$

When the response set $\mathcal{S}(x,\varepsilon)$ is a lower semi-continuous multivalued function, and the optimal value function $\varphi(x)$ is continuous, the proposed continuous functions $\beta(\cdot)$ in $\mathcal{X}$ ensure the existence of a solution to the {\SWDBO}~\citep{jia2011new}. These properties are satisfied for the linear follower problems with bounded feasible sets used in our experiments.
\begin{proposition}\label{prop:swd}
    The optimal leader value $\Theta_{\I}(\mu_{\beta_{\text{Prp}}(\cdot)})$ is a lower-bound to $\Theta_{\I}(\mu_{\beta(\cdot)})$ for any convex $\beta(\cdot)$ (e.g., $\beta_{\text{Frg}}(\cdot)$ and $\beta_{\text{Frg-p}}(\cdot)$) and an upper-bound to $\Theta_{\I}(\mu_{\beta(\cdot)})$ for any concave $\beta(\cdot)$ (e.g., $\beta_{\text{Str}}(\cdot)$ and $\beta_{\text{Str-p}}(\cdot)$).
\end{proposition}

\subsection{Generalized Decision-Dependent Measures}\label{subsec:modeling-ibl}

Beyond the neutral measure, we define endogenous measures over the full support $\mathcal{S}(x,\varepsilon)$. To enable scale-invariant comparisons, we define the normalized leader objective:
\begin{equation}\label{eq:normalized-leader-value}
    \Tilde{F}(x,\hat{y}) = \begin{cases} 
    \frac{F(x,\hat{y}) - \tfrac{1}{2} [F_{\W}(x,\varepsilon) + F_{\Snew}(x,\varepsilon)]}{F_{\W}(x,\varepsilon) - F_{\Snew}(x,\varepsilon)} \in  \left[-\tfrac{1}{2}, \tfrac{1}{2}\right]  & \text{if } F_{\Snew}(x,\varepsilon) < F_{\W}(x,\varepsilon), \\
    0 & \text{if } F_{\Snew}(x,\varepsilon) = F_{\W}(x,\varepsilon),
    \end{cases} \quad \forall\ x \in \mathcal{X},\ \hat{y} \in \mathcal{S}(x,\varepsilon).
\end{equation}
This function centers outcomes relative to the midpoint of the optimistic and pessimistic values. If $F_{\Snew}(x,\varepsilon) = F_{\W}(x,\varepsilon)$, all near-optimal solutions yield the same leader value and thus have equal probability.

For each $x \in \mathcal{X}$, we define the probability measure $\mu_x^{\pi(\cdot)}$ over $\mathcal{S}(x,\varepsilon)$ proportional to a nonnegative, measurable function of both the leader’s objective value and the follower’s optimal value, $\pi \left( \Tilde{F}(x,\hat{y}) ,\ \tilde{\varphi}(x) \right)$, such that ${\pi: \left[-\tfrac{1}{2}, \tfrac{1}{2}\right] \times [0,1] \to [\pi_{min},\pi_{max}]}$ with ${0 < \pi_{min} < \pi_{max} < \infty}$. The endogenous probability measure is thus obtained by normalization:
\begin{equation}\label{eq:measure-pi}
    \mu^{\pi(\cdot)}_x(A) = \frac{\int_{\hat{y} \in A}  \pi \left( \tilde{F}(x,\hat{y}) , \tilde{\varphi}(x) \right) \, d \hat{y} }{\int_{\hat{y} \in \mathcal{S}(x,\varepsilon)} \pi \left( \tilde{F}(x,\hat{y}) , \tilde{\varphi}(x) \right) \, d \hat{y} }, \qquad A \subseteq \mathcal{S}(x,\varepsilon),
\end{equation}
with $\mu^{\pi(\cdot)} = \cup_{x \in \mathcal{X}} \left\{ \mu_x^{\pi(\cdot)} \right\}$. If $\pi(\cdot)$ is constant (e.g., ${\pi_{\text{Ntr}}\! \left( \Tilde{F}, \tilde{\varphi} \right) = 1}$), this reduces to the neutral measure $\mu_x^{\text{Ntr}}$. Otherwise, $\pi(\cdot)$ dictates behavior as described previously: decreasing in $F(x,\hat{y})$ for cooperative settings ($\tilde\varphi(x)$ near $0$) and increasing for adversarial ones ($\tilde\varphi(x)$ near $1$).

\begin{figure}[t]
    \centering
    \caption{Proposed cooperation functions $\pi\!\left(\tilde{F}(x,\hat{y}), \tilde{\varphi}(x)\right)$ for $\tilde{F}(x,\hat{y}) \in \left[-\tfrac{1}{2}, \tfrac{1}{2}\right]$ and fixed values $\tilde{\varphi}(x) = 0.2$ (most cooperative), $0.4$, $0.6$, and $0.8$ (least cooperative)}. \label{fig:pi-functions}
    %
    \pgfplotsset{
        Vone/.style   ={thick, solid},
        Vtwo/.style   ={thick, dashed},
        Vthree/.style ={thick, dotted},
        Vfour/.style  ={thick, dash dot}
    }
    %
    \begin{subfigure}{0.22\textwidth}
        \centering
        \begin{tikzpicture}
        \begin{axis}[
            xlabel={$\Tilde{F}(x,\hat{y})$},
            xmin=-0.55,xmax=0.55,
            samples=100,
            xtick={-0.5,0.5}, xticklabels={$-\tfrac{1}{2}$,$\tfrac{1}{2}$},
            ytick={0,2.25,4.5}, yticklabels={$\pi_{\min}$,$\pi_c$,$\pi_{\max}$},
            height=3.2cm,
            grid=both
        ]
        \def\ymid{2.25}
        \addplot[Vone,domain=-0.5:0.5] {\ymid - 4.5*x};
        \addplot[Vtwo,domain=-0.5:0.5] {\ymid - 1.5*x};
        \addplot[Vthree,domain=-0.5:0.5] {\ymid + 1.5*x};
        \addplot[Vfour,domain=-0.5:0.5] {\ymid + 4.5*x};
        \end{axis}
        \end{tikzpicture}
        \caption{Proportional $\pi_{\text{Prp}}(\cdot)$}
    \end{subfigure}
    \hfill
    %
    \begin{subfigure}{0.22\textwidth}
        \centering
        \begin{tikzpicture}
        \begin{axis}[
            xlabel={$\Tilde{F}(x,\hat{y})$},
            xmin=-0.55,xmax=0.55,
            samples=100,
            xtick={-0.5,0.5}, xticklabels={$-\tfrac{1}{2}$,$\tfrac{1}{2}$},
            ytick={0.0,1.0}, yticklabels={$\pi_{\min}$,$\pi_{\max}$},
            height=3.2cm,
            grid=both
        ]
        \def\ymid{1}
        \addplot[Vone,domain=-0.5:0.5] {\ymid*exp(-4.5*x -4.5*0.5)};
        \addplot[Vtwo,domain=-0.5:0.5] {\ymid*exp(-1.5*x -1.5*0.5)};
        \addplot[Vthree,domain=-0.5:0.5] {\ymid*exp(1.5*x -1.5*0.5)};
        \addplot[Vfour,domain=-0.5:0.5] {\ymid*exp(4.5*x -4.5*0.5)};
        \end{axis}
        \end{tikzpicture}
        \caption{Mod. Fragile $\pi_{\text{Frg}}(\cdot)$}
    \end{subfigure}
    \hfill
    %
    \begin{subfigure}{0.22\textwidth}
        \centering
        \begin{tikzpicture}
        \begin{axis}[
            xlabel={$\Tilde{F}(x,\hat{y})$},
            xmin=-0.55,xmax=0.55,
            samples=100,
            xtick={-0.5,0.5}, xticklabels={$-\tfrac{1}{2}$,$\tfrac{1}{2}$},
            ytick={0,1}, yticklabels={$\pi_{\min}$,$\pi_{\max}$},
            height=3.2cm,
            grid=both
        ]
        \def\xmid{0.5}
        \addplot[Vone,domain=-0.5:0.5] {pow(\xmid - x, 1 + 4.5)};
        \addplot[Vtwo,domain=-0.5:0.5] {pow(\xmid - x, 1 + 1.5)};
        \addplot[Vthree,domain=-0.5:0.5] {pow(\xmid + x, 1 + 1.5)};
        \addplot[Vfour,domain=-0.5:0.5] {pow(\xmid + x, 1 + 4.5)};
        \end{axis}
        \end{tikzpicture}
        \caption{Pron. Fragile $\pi_{\text{Frg}^{+}} ( \cdot )$}
    \end{subfigure}
    \hfill
    %
    \begin{subfigure}{0.22\textwidth}
        \centering
        \begin{tikzpicture}
        \begin{axis}[
            xlabel={$\Tilde{F}(x,\hat{y})$},
            xmin=-0.55,xmax=0.55,
            samples=200,
            xtick={-0.5,0.5}, xticklabels={$-\tfrac{1}{2}$,$\tfrac{1}{2}$},
            ytick={0,1}, yticklabels={$\pi_{\min}$,$\pi_{\max}$},
            height=3.2cm,
            grid=both,
            legend style={
            at={(0.5,1.05)},
            anchor=south,
            draw=none,
            /tikz/every even column/.append style={column sep=0.8cm}
        }
        ]
        \addplot[Vone,domain=-0.5:0.5] {pow(0.5 - x,0.25)};
        \addplot[Vtwo,domain=-0.5:0.5] {pow(0.5 - x,0.084)};
        \addplot[Vthree,domain=-0.5:0.5] {pow(0.5 + x,0.084)};
        \addplot[Vfour,domain=-0.5:0.5] {pow(0.5 + x,0.25)};
        \end{axis}
        \end{tikzpicture}
        \caption{Sturdy $\pi_{\text{Str}}(\cdot)$}
    \end{subfigure}
\vspace{0.1cm}
    \begin{tikzpicture}
    \begin{axis}[
        hide axis,
        xmin=0,xmax=1,ymin=0,ymax=1,
        legend columns=5,
        legend style={
            at={(0.5,1.05)},
            anchor=south,
            draw=none
        }
    ]
    \addlegendimage{Vone}\addlegendentry{ $\tilde{\varphi}(x) = 0.2 $ }
    \addlegendimage{Vtwo}\addlegendentry{ $\tilde{\varphi}(x) = 0.4 $ }
    \addlegendimage{Vthree}\addlegendentry{ $\tilde{\varphi}(x) = 0.6 $ }
    \addlegendimage{Vfour}\addlegendentry{ $\tilde{\varphi}(x) = 0.8 $ }
    \end{axis}
    \end{tikzpicture}
\end{figure}
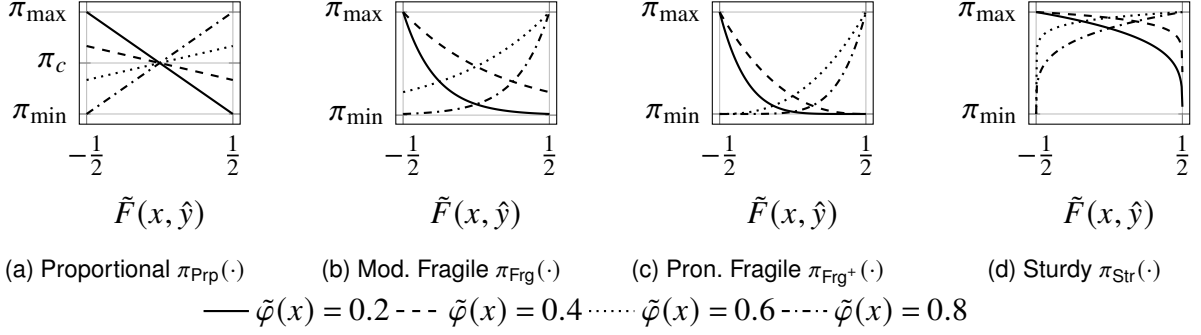

We formalize three generalized forms of follower cooperation: Proportional, Fragile (Moderate and Pronounced or Polynomial), and Sturdy. The mathematical definitions for these generalized measures--$\pi_{\text{Prp}}(\cdot)$, $\pi_{\text{Frg}}(\cdot)$, $\pi_{\text{Frg}^{+}}(\cdot)$, and $\pi_{\text{Str}}(\cdot)$--are presented in Table~\ref{tab:cooperation_functions} alongside their dichotomous counterparts. Here, $\gamma \geq 0$ scales the growth rate for fragile and sturdy behaviors, ${m(\tilde{\varphi}(x)) = 2\tilde{\varphi}(x) - 1 \in [-1,1]}$ is the continuous cooperation factor, and ${s(\tilde{\varphi}(x)) = \operatorname{sgn}(m(\tilde{\varphi}(x)))}$. The constant $\pi_c$ (e.g., $1/2 + \epsilon$) ensures $\pi_{\text{Prp}}(\cdot) > 0$. We denote the set of generalized functions as $\Pi = \{\pi_{\text{Ntr}}(\cdot),\pi_{\text{Prp}}(\cdot),\pi_{\text{Frg}}(\cdot),\pi_{\text{Frg}^{+}}(\cdot),\pi_{\text{Str}}(\cdot)\}$. All generalized cooperation types are illustrated in Figure~\ref{fig:pi-functions}. In Example~\ref{ex:behavior-ibl} (Appendix~\ref{app:examples}), we illustrate the impact of this type of behavior on the optimal decisions and values for both players. The example shows that even Generalized Proportional behavior can lead to different solutions compared to the neutral measure.
\begin{table}[ht]
\TABLE{Definitions of dichotomous ($\beta$) and generalized ($\pi$) follower cooperation behaviors.\label{tab:cooperation_functions}}
{
\begin{tabular}{lllc}
\hline
\textbf{Behavior} & \textbf{Dichotomous Measure $\beta(\tilde{\varphi}(x))$} & \textbf{Generalized Measure $\pi(\tilde{F}(x,\hat{y}),\tilde{\varphi}(x))$} & \textbf{Parameters} \\ \hline
\text{Proportional} & $1 - \tilde{\varphi}(x)$ & $\pi_c + m(\tilde{\varphi}(x)) \Tilde{F}(x,\hat{y})$ & --- \\
\text{Threshold} & $\frac{1}{1 + \exp(\delta[\tilde{\varphi}(x) - 1/2])}$ & --- & $\delta > 0$ \\
\text{Sturdy} & $\frac{1 - \exp(-a[1 - \tilde{\varphi}(x)])}{1 - \exp(-a)}$ & $\left( \tfrac{1}{2} + s(\tilde{\varphi}(x)) \Tilde{F}(x,\hat{y}) \right)^{\frac{1}{\gamma} |m(\tilde{\varphi}(x))|}$ & $a > 0, \gamma \ge 1$ \\
\text{Fragile} & $1 - \frac{1 - \exp(-a\tilde{\varphi}(x))}{1 - \exp(-a)}$ & $\exp\left(\gamma m(\tilde{\varphi}(x)) \Tilde{F}(x,\hat{y}) - \frac{\gamma}{2} |m(\tilde{\varphi}(x))|\right)$ & $a > 0, \gamma > 0$ \\
\text{Polynom. Sturdy} & $1 - \tilde{\varphi}(x)^p$ & --- & $p > 1$ \\
\text{Polynom. Fragile} & $(1 - \tilde{\varphi}(x))^p$ & $\left( \tfrac{1}{2} + s(\tilde{\varphi}(x)) \Tilde{F}(x,\hat{y}) \right)^{1 + \gamma |m(\tilde{\varphi}(x))|}$ & $p > 1, \gamma > 0$ \\ \hline
\end{tabular}
}
{}
\end{table}

Similarly to the dichotomous case, suppose that the response set $\mathcal{S}(x,\varepsilon)$ is a lower semi-continuous multivalued function, and the optimal value function $\varphi(x)$ is continuous. Then, if $\pi\!(\cdot)$ is continuous on $\mathcal{S}(x,\varepsilon)$ (as in the proposed behaviors), and the endogenous measure $\mu_x^{\pi(\cdot)}$ is absolutely continuous with respect to the Lebesgue measure, then its density depends continuously on $\mathcal{X}$, thereby guaranteeing the existence of a solution for the {\IBO}~\citep{lina1996hierarchical}.

\section{Methodology}\label{sec:methodology}

The term ${\mathbb{E}_{\boldsymbol{\hat{y}}\!\sim\!\mu_x}\!\left[F(x,\boldsymbol{\hat{y}})\right]}$ generally lacks a closed-form expression. This makes {\DE} formulations for {\IBO} (Program~\eqref{Prog:intermediate}), even when available, highly nonlinear and nonconvex, and direct sampling nontrivial. To handle this, we adopt the random variable transformation paradigm of~\citet{bazotte2025solving} to map the stochastic program with endogenous uncertainty into an equivalent program with exogenous uncertainty. The resultant Transformed {\IBO} ({\TIBO}) is defined as:
\begin{equation}\label{Prog:transf-intermediate}
    \Theta_{\I} \! \left( {\mu} \right) = \min_{x \in \mathcal{X}} \left\{ H(x) + g^{\mu} (x) \right\}, \text{ with } g^{\mu} (x)  =  \mathbb{E}_{\boldsymbol{\zeta}\sim \rho \left[ F \left( x, \boldsymbol{\hat{y}} =  t_{\mu} (x, \boldsymbol{\zeta}) \right) \right]} ,
\end{equation}
where $\boldsymbol{\zeta}$ is a continuous exogenous random vector with realizations ${\zeta \in \mathcal{Z}}$ and density $\varrho$, and ${\boldsymbol{\hat{y}} = t_{\mu} (x, \boldsymbol{\zeta})}$ is an appropriate transformation (specified in the following subsections) mapping the leader's decisions and the exogenous uncertainty to the endogenous response $\boldsymbol{\hat{y}}\sim \mu_x$ with support $\mathcal{S}(x,\varepsilon)$. We similarly define the expected follower objective $g_f^\mu(x) = \mathbb{E}_{\boldsymbol{\zeta} \sim \varrho} \left[ f\!\left(x,t_{\mu}(x,\boldsymbol{\zeta})\right) \right]$.

Enumerating all realizations of the continuous random vector $\boldsymbol{\zeta}$ is infeasible; however, its exogeneity enables the application of sampling-based methods such as the {\SAA} to approximate the expectation in Program~\eqref{Prog:transf-intermediate}. The {\SAA} program of the {\TIBO} ({\SAA} {\TIBO}) is then defined as:
\begin{equation}\label{Prog:saa-transf-intermediate}
    \overline{\Theta}_{\I}^{N} \! (\mu) = \min_{x \in \mathcal{X}} \left\{ H(x) + \hat{g}^{\mu}_N (x) \right\}, \, \text{with} \, \hat{g}^{\mu}_N (x)  = \! \frac{1}{N} \!\sum_{\omega \in \mathcal{W}_N} \! F \! \left( x, \hat{y}_\omega = t_{\mu} (x, \zeta_\omega) \right),
\end{equation}
where $\hat{g}^{\mu}_N (x)$ is the {\SAA} estimator of $g^{\mu} (x)$, with $\mathcal{W}_N = \{\omega_1,\ldots,\omega_N\}$ a set of $N$ independent and identically distributed samples of $\boldsymbol{\zeta}$. We also define ${\hat{g}^{\mu}_{f,N} (x) \! = \! 1/N \! \cdot \! \sum_{w \mathcal{W}_N} \! f\!\left(x,\hat{y}_\omega = t_{\mu}(x,\zeta_\omega) \right)}$, the unbiased estimator of $g^\mu_f(x)$. We denote by $\hat{y}_\omega \in \mathcal{S}(x,\varepsilon)$ the realization of $\boldsymbol{\hat{y}}$ corresponding to the sample $\zeta_\omega$ through the transformation $t_{\mu} (\cdot)$. The {\SAA} method proceeds in two steps~\citep[see chap. 5]{shapiro2021lectures}. First, we solve {\SAA}~Program~\eqref{Prog:saa-transf-intermediate} to obtain a feasible solution $\overline{x}^{\mu}_N$ with objective value $\overline{\Theta}_{\I}^N\!(\mu)$. Second, we evaluate the upper bound $g^{\mu}\left(\overline{x}_N^\mu\right)$ either exactly (if the support is enumerable) or via a large independent sample $N' \gg N$ to compute $\hat{g}^{\mu}_{N'}(\overline{x}^{\mu}_N)$. We note that $\overline{\Theta}_{\I}^N\!(\mu)$ and $\hat{g}^{\mu}_{N'}\! \left(\overline{x}^{\mu}_N\right)$ (or $g^{\mu}\! \left(\overline{x}^{\mu}_N\right)$) are lower and and upper bound estimators of the true optimal value $\Theta_{\I}\!\left( \mu \right)$. While solving multiple {\SAA} instances with distinct sets $\mathcal{W}_N$ can improve the lower bound estimator, doing so is expected to be computationally prohibitive given the hardness of bilevel programs; hence, we solve a single instance.

Developing these transformations is challenging because the support $\mathcal{S}(x,\varepsilon)$ is implicitly defined by the follower's optimization problem. To our knowledge, we are the first to extend this transformation framework to bilevel models with decision-dependent tie-breaking. Next, we derive $\boldsymbol{\zeta}$ and $t_{\mu}(x,\boldsymbol{\zeta})$ for the dichotomous measures (Section~\ref{subsec:methodology-swd}) and generalized measures (Section~\ref{subsec:methodology-ibl}).

\subsection{Strong Weak Decision-Dependent Intermediate} \label{subsec:methodology-swd}

For the {\SWDBO}, we apply the inverse transformation technique of~\citet{bazotte2025solving} to construct the decision-dependent Bernoulli random variable in Program~\eqref{Prog:swd}. Let ${\boldsymbol{\zeta} \sim U(0,1)}$ with ${\mathcal{Z} = [0,1]}$. For any realization ${\zeta \in \mathcal{Z}}$, the inverse transformation ${t_{\mu_{\beta(\cdot)}}(x,\zeta)}$ returns ${\hat{y}^{\Snew}(x,\varepsilon)}$ if ${\zeta \le \beta\!\left(\tilde{\varphi}(x)\right)}$, and ${\hat{y}^{\W}(x,\varepsilon)}$ otherwise, where $\hat{y}^{\Snew}(x,\varepsilon)$ and $\hat{y}^{\W}(x,\varepsilon)$ denote the optimistic and pessimistic solutions, respectively (see Section~\ref{subsec:literature-background}). The strict monotonicity of $\beta \! \left(\tilde{\varphi}(x) \right)$ (Assumption~\ref{assump:mono}) guarantees the existence of a unique inverse $\beta^{-1} ( \zeta )$ for $\zeta \in \mathcal{Z}$, which allows the transformation to be equivalently defined through this function. In particular, the mapping ${t_{\mu_{\beta(\cdot)}}(x,\zeta)}$ returns ${\hat{y}^{\Snew}(x,\varepsilon)}$ if ${\beta^{-1} \!(\zeta) \ge \tilde{\varphi}(x)}$, and ${\hat{y}^{\W}(x,\varepsilon)}$ otherwise.

To specify the {\SAA} formulation of the transformed program, we model the functions $\hat{y}^{\Snew}(x,\varepsilon)$ and $\hat{y}^{\W}(x,\varepsilon)$ and the proposed transformation. For that, we introduce variables ${y_{\Snew} \in \mathbb{R}^{n_y}_+}$ and ${y_{\W} \in \mathbb{R}_+^{n_y}}$ to represent ${\hat{y}^{\Snew}(x,\varepsilon)}$ and ${\hat{y}^{\W}(x,\varepsilon)}$, respectively, and, similarly, define variables $F_{\Snew}$ and $F_{\W}$ representing $F_{\Snew}(x,\varepsilon)$ and $F_{\W}(x,\varepsilon)$. To correctly define these variables, and also the normalized follower optimal value $\tilde{\varphi}(x)$, we define variables $y_{\texttt{O}}$ to be an optimal solution of the follower's problem. We assume the scenarios are ordered in decreasing values of $\beta^{-1}(\cdot)$; that is, ${\beta^{-1}(\zeta_{\omega_i}) \ge \beta^{-1}(\zeta_{\omega_{i+1}})}$, ${i =1,\ldots,N-1}$. We then introduce the variables ${z_{\omega_i} \in \{0,1\}}$ for $i =1,\ldots,N$ to identify the threshold where ${\tilde{\varphi}(x)}$ crosses the ordered sequence $\left\{\beta^{-1}(\zeta_{\omega_i})\right\}_{i=1}^N$. Specifically, $z_{\omega_i} = 1$ if ${\beta^{-1}(\zeta_{\omega_i}) \ge \tilde{\varphi}(x)}$ and ${\beta^{-1}(\zeta_{\omega_{i+1}}) < \tilde{\varphi}(x)}$, and $z_{\omega_i} = 0$ otherwise. The {\SAA} program of the {\TIBO} for dichotomous decision-dependent measures is:
\begin{subequations}\label{Prog:saa-transformed-swd}
    \begin{align}
        \min_{x \in \mathcal{X}, \, z \in \{0,1\}^{N}} \mbox{: } & H(x) + F_{\W} + \sum_{\omega_i \in \mathcal{W}_N} \! \tfrac{i}{N} \left(F_{\Snew} -  F_{\W} \right) z_{\omega_i} \\
        \mbox{s.t.: } \label{constr:def-follower} & y_{\texttt{O}} \in \mathcal{S}(x,0), \quad \tilde{\varphi}(x) = \tfrac{f(x,y_{\texttt{O}}) - \varphi_{lb}}{\varphi_{ub} - \varphi_{lb}}, \\
        \label{constr:def-opt} & F_{\Snew} = F(x,y_{\Snew}), \quad y_{\Snew} \in \left\{ y \in \mathcal{Y}(x) : f(x,y) \le f(x,y_{\texttt{O}}) + \varepsilon  \right\}, \\
        \label{constr:def-pes} & F_{\W} = F(x,y_{\W}), \quad y_{\W} \in  \arg\max_y \left\{ F(x,y): f(x,y) \leq f(x,y_{\texttt{O}}) + \varepsilon, y \in \mathcal{Y}(x) \right\}, \\
        \label{constr:transf}& \sum_{\omega \in \mathcal{W}_N} z_{\omega} \le 1, \quad \tilde{\varphi}(x) \le \beta^{-1}\!\left(\zeta_{\omega} \right) + \left( 1 - \beta^{-1}\!\left( \zeta_{\omega} \right) \right) \left(1 - z_{\omega} \right) \; \forall\ \omega \in \mathcal{W}_N .
    \end{align}
\end{subequations}
Constraint~\eqref{constr:def-follower} 
guarantees that $y_{\texttt{O}}$ is an optimal follower response. Constraints~\eqref{constr:def-opt} and~\eqref{constr:def-pes} define the $\varepsilon$-near-optimal optimistic and pessimistic follower solutions, respectively. In the special case where $\varepsilon=0$, the variables $y_{\Snew}$ can be directly constrained to belong to $\mathcal{S}(x,\varepsilon)$, and the variables $y_{\texttt{O}}$ are no longer required. The pessimistic solution in Constraint~\eqref{constr:def-pes} is modeled using a modified version of the approach of~\citet{zeng2020practical}. We thus have the following proposition:
\begin{proposition}\label{prop:def-opt-pes}
    Given $y_{\texttt{O}} \in \mathcal{S}(x,0)$, Constraints~\eqref{constr:def-opt} and~\eqref{constr:def-pes} correctly determine ${F_{\Snew} = \min_{\hat{y} \in \mathcal{S}(x,\varepsilon)} F(x,\hat{y})}$ and ${F_{\W} =\max_{\hat{y} \in \mathcal{S} (x,\varepsilon)} F(x,\hat{y})}$.
\end{proposition}
As a result, the {\SAA} Program~\eqref{Prog:saa-transformed-swd} exactly identifies the follower’s $\varepsilon$-near-optimal optimistic and pessimistic values, without the correction scheme in~\citet{zeng2020practical}. Constraints~\eqref{constr:transf}, together with the objective function, define the variables $z_\omega$ and model the transformation $t_{\mu_{\beta(\cdot)}}(\cdot)$. The objective minimizes the expected leader value by accounting for the probability of selecting either the optimistic or pessimistic solution, leading to the bilinear terms $F_{\Snew} z_{\omega}$ and $F_{\W} z_{\omega}$, which can be exactly linearized using McCormick envelopes~\citep{mccormick1976computability}. The formulation handles any function satisfying Assumption~\ref{assump:mono}. Standard reformulations apply to the nested follower problems; for linear followers, {\KKT} conditions yield complementarity constraints, producing a single-level {\MILP} solvable by standard solvers regardless of the shape of $\beta(\tilde{\varphi}(x))$.

As a baseline, and using Proposition~\ref{prop:def-opt-pes}, the {\DE} of the {\SWDBO} ({\DE} {\SWDBO}) is:
\begin{equation}\label{Prog:de-swd}
    \min_{x \in \mathcal{X}} \left\{ H(x) + \beta \left( \tilde{\varphi}(x) \right)  F_{\Snew}  + \left[1 - \beta \left( \tilde{\varphi}(x) \right) \right] F_{\W}: y_{\texttt{O}} \in \mathcal{S}(x,0), \eqref{constr:def-opt}, \eqref{constr:def-pes} \right\}.
\end{equation}
Unlike our {\MILP} {\SAA} program, this {\DE} formulation remains a nonlinear, nonconvex bilevel program due to the form of $\beta \left(\tilde{\varphi}(x) \right)$. This program generalizes the {\DE} of the {\SWBO} ({\DE} {\SWDBO}), in which case the cooperation function $\beta \left(\tilde{\varphi}(x) \right)$ is instead a fixed value $\beta \in [0,1]$.

\subsection{Generalized Decision-Dependent Intermediate} \label{subsec:methodology-ibl}

For the generalized measures in Section~\ref{subsec:modeling-ibl}, the expectation $\mathbb{E}_{\boldsymbol{\hat{y}}\sim\mu_x}[ F(x,\boldsymbol{\hat{y}}) ]$ typically does not admit a closed-form expression, making explicit {\DE} formulations unavailable and computationally challenging. Approximating this expectation requires sampling from $\mathcal{S}(x,\varepsilon)$, which involves constrained procedures and generally relies on specialized methods such as rejection sampling or {\MCMC}. To address this challenge, we extend the Integer {\sf L}-shaped method~\citep{laporte1993integer} to solve the {\SAA} {\TIBO} when the leader decisions are binary, i.e., $\mathcal{X} \subseteq \{0,1\}^{n_x}$. This implicitly embeds the sampling-based transformation within a branch-and-cut scheme. As in the previous case, the variables $y_{\texttt{O}} \in \mathbb{R}_+^{n_y}$ represent an optimal follower solution. We introduce variables $y_{\texttt{NO}} \in \mathbb{R}_+^{n_y}$ to denote an $\varepsilon$-near-optimal follower solution, and $\theta$ to capture the estimator value $\hat{g}^{\mu}_N(x)$. The resulting {\SAA} formulation of the {\TIBO} for generalized measures is given as follows:
\begin{subequations}\label{Prog:saa-transf-int-cuts}
    \begin{align}
       \min_{x \in \mathcal{X} \subseteq \{0,1\}^{n_x}} \mbox{: } & H(x) + \theta \\
       \mbox{s.t.: } \label{constr:def-follower-nophi} & y_{\texttt{O}} \in \mathcal{S}(x,0),  \\
       \label{constr:lb-opt}  & y_{\texttt{NO}} \in \left\{ y \in \mathcal{Y}(x) : f(x,y) \le f(x,y_{\texttt{O}}) + \varepsilon  \right\}, \quad \theta \ge F(x,y_{\texttt{NO}}), \\
       \label{constr:integer-cut} & \theta \ge \left(\hat{g}^{\mu}_N(\tilde{x}) - F_{lb} \right) \left( \sum_{j \in \mathcal{J}(\tilde{x})} x_j - \sum_{j \notin \mathcal{J}(\tilde{x})} x_j - \left| \mathcal{J}(\tilde{x}) \right| \right) + \hat{g}^{\mu}_N\!(\tilde{x}) \; \forall\ \tilde{x} \in \mathcal{X}^r \subseteq \mathcal{X},     
    \end{align}
\end{subequations}
where $F_{lb}$ is a lower bound on $F(x,\hat{y})$ for all $x \in \mathcal{X}$ and $\hat{y} \in \mathcal{S}(x,\varepsilon)$, i.e., ${F_{lb}=\min_{x \in \mathcal{X}} \{F_{\Snew}: \eqref{constr:def-follower-nophi},\eqref{constr:def-opt} \}}$, and $\mathcal{J}(\tilde{x}) = \{j : \tilde{x}_j = 1\}$. Constraint~\eqref{constr:def-follower-nophi} ensures that $y_{\texttt{O}}$ is an optimal follower solution. Constraints~\eqref{constr:lb-opt} define the variables $y_{\texttt{NO}}$, whose associated leader objective value $F(x,y_{\texttt{NO}})$ provides a valid lower bound on the auxiliary variable $\theta$. The cuts~\eqref{constr:integer-cut} are dynamically generated within a branch-and-cut scheme. At iteration $r$, $\mathcal{X}^r \subseteq \mathcal{X}$ denotes the set of feasible leader solutions for which cuts have been added, with $\mathcal{X}^0 = \varnothing$. Although feasible $y_{\texttt{NO}}$ may yield values of $F(x,y_{\texttt{NO}})$ that are either above or below the estimator $\hat{g}^{\mu}_N(x)$, the epigraph formulation with $\theta$, combined with cuts~\eqref{constr:integer-cut}, ensures that $\theta$ correctly approximates the estimator. Finally, the formulation extends to bounded integer leader variables via a standard binary expansion.

The computation of $\hat{g}^{\mu}_N(x)$ for a given $x$ depends on $\mathcal{S}(x,\varepsilon)$. We consider linear follower problems in both $x$ and $y$, defined by $f_{\text{Lin}}(x,y) = e^\top y$, $\mathcal{Y}_{\text{Lin}}(x) = \left\{ y \ge 0 : Dy \le a - Cx \right\}$, and ${\varphi_{\text{Lin}}(x) = \min_y \left\{ e^\top y : Dy \le a - Cx, \, y \ge 0 \right\}}$, where $D \in \mathbb{R}^{m_2 \times n_y}$, $a \in \mathbb{R}^{m_2}$, $C \in \mathbb{R}^{m_2 \times n_x}$, and $e \in \mathbb{R}^{n_y}$. The $\varepsilon$-optimal response set is ${\mathcal{S}_{\text{Lin}}(x,\varepsilon) \!=\! \left\{ \hat{y} \ge 0 : e^\top \hat{y} \le \varphi_{\text{Lin}}(x) + \varepsilon,  D\hat{y} \le a - Cx \right\}}$, a bounded polyhedral set for fixed $x$ and compact $\mathcal{Y}_{\text{Lin}}(x)$. The {\KKT} conditions can be used to model Constraints~\eqref{constr:def-follower-nophi}, yielding a single-level reformulation of {\SAA}~Program~\eqref{Prog:saa-transf-int-cuts}. 

Algorithm~\ref{alg:cut-separation} presents the separation scheme for Cuts~\eqref{constr:integer-cut} for this class of linear lower-level problems, where $\boldsymbol{1}$ denotes the all-ones vector in $\mathbb{R}^{n_y}$ and $\|D\|_{\mathrm{row},2}$ the vector of Euclidean norms of the rows of $D$. An integer {\sf L}-shaped cut can be separated for any leader decision $x \in \mathcal{X}$ such that $\mathcal{S}(x,\varepsilon) \neq \varnothing$. Accordingly, within the branch-and-bound tree, the separation algorithm is invoked whenever a solution $(x, y_{\texttt{O}})$ with $x \in \mathcal{X}$ and $y_{\texttt{O}} \in \mathcal{Y}_{\mathrm{Lin}}(x)$ (i.e., satisfying primal feasibility) is found, even if $y_{\texttt{O}} \notin \mathcal{S}(x,0)$ (i.e., the {\KKT} complementarity conditions are not satisfied). This ensures that $\mathcal{S}(x,\varepsilon)$ is nonempty when the separation algorithm is executed. Consequently, the optimal follower value is computed accordingly, depending on whether $y_{\texttt{O}} \in \mathcal{S}(x,0)$ (Step~1 of the algorithm).

\begin{algorithm}[t]
    \small
    \caption{Separation Scheme for Cuts~\eqref{constr:integer-cut}} \label{alg:cut-separation} 
    \small 
    \Input{$\varepsilon$, $\varphi_{lb}$, $\varphi_{ub}$, $\gamma$, $K$, $\pi(\cdot) \in \Pi$, $\mu=\mu_{\pi(\cdot)}$, $\mathcal{X}^{r-1}$, $x \in \mathcal{X}$, $y_{\texttt{O}} \in \mathcal{Y}_{\mathrm{Lin}}(x)$;}
    \Output{$\mathcal{X}^r$, $\hat{g}^{\mu}_N\!(x)$;} 
    
    \Stepone{If $y_{\texttt{O}} \notin \mathcal{S}_{\mathrm{Lin}}(x,0)$, solve $f_{\texttt{O}} = \varphi_{\mathrm{Lin}}(x)$; otherwise, set $f_{\texttt{O}} = f_{\mathrm{Lin}}(x, y_{\texttt{O}})$. Set $\tilde{\varphi}(x) = \tfrac{f_{\texttt{O}} - \varphi_{lb}}{\varphi_{ub} - \varphi_{lb}}$;} 

    \Steptwo{Compute the $\varepsilon$-near optimistic $F_{\Snew}(x,\varepsilon)$ and pessimistic $F_{\W}(x,\varepsilon)$ leader values for $x$ to define $\tilde{F}(x,\hat{y})$ as in Equation~\eqref{eq:normalized-leader-value};} 

    \Stepthree{Compute the Chebyshev center of $\mathcal{S}_{\mathrm{Lin}}(x,\varepsilon)$: 
    \[
        (\hat{y}_{\texttt{C}}, k_{\texttt{C}}) \in \arg\max_{y,k} \left\{ k : 
        Dy + \|D\|_{\text{row},2} \, k \leq a - Cx, \;
        e^\top y + \|e\|_2 \, k \leq f_{\texttt{O}} + \varepsilon, \;
        y - \boldsymbol{1} k \geq 0 
        \right\};
    \]}

    \Stepfour{Compute $\hat{g}^{\mu}_N(x) = \tfrac{1}{N} \sum_{\omega \in \mathcal{W}_N} F(x,\hat{y}_\omega)$ with samples $\{\hat{y}_{\omega}\}_{\omega \in \mathcal{W}_N}$ drawn from $\mathcal{S}_{\mathrm{Lin}}(x,\varepsilon)$ proportionally to $\pi \left( \tilde{F}(x,\hat{y}), \tilde{\varphi}(x) \right)$ with parameter $\gamma$ via the {\MCMC} Hit-and-Run, initialized at $\hat{y}_{\texttt{C}}$;}

    \Stepfive{Define set $\mathcal{J}(x)$ and Cut~\eqref{constr:integer-cut} for fixed $x$, and set $\mathcal{X}^{r} \gets \mathcal{X}^{r-1} \cup \{ x \}$.}
\end{algorithm}

In Step~4 of Algorithm~\ref{alg:cut-separation}, we compute $\hat{g}^{\mu}_N(x)$ using samples from $\mathcal{S}_{\text{Lin}}(x,\varepsilon)$ generated via {\MCMC} Hit-and-Run, which is well-suited for sampling from polyhedral sets~\citep[see, e.g.,][]{lovasz2006hit} (see Appendix~\ref{app:details-har} for details on this procedure). More generally, Hit-and-Run applies to convex sets~\citep{belisle1993hit}, and the proposed branch-and-cut framework extends to broader uncertainty sets $\mathcal{S}(x,\varepsilon)$ by combining an appropriate sampling method with the cut separation step.

The transformation is implicitly defined by the Markov chain induced by the Hit-and-Run dynamics. While the endogenous samples $\{\hat{y}_{\omega}\}_{\omega \in \mathcal{W}_N}$ are dependent, the underlying exogenous random vectors guiding the steps are independently distributed, ensuring convergence and consistent estimation. This allows the exact branch-and-cut procedure to optimally solve the {\SAA} Program~\eqref{Prog:saa-transf-int-cuts}.

\section{Experimental Results}\label{sec:experiments}

The goals of our experiments are twofold: (i) to assess the computational performance of our {\SAA}-based methods for the {\IBO}; (ii) to examine how the leader’s and follower’s optimal values vary under different follower behaviors, including the impact of follower behavior misspecification.

We adapt 170 instances of the BOBILib library~\citep{BOBILib:2026} (110 instances from the \texttt{denegre} set~\citep{denegre2011interdiction} and 60 from the \texttt{xuwang} set~\citep{xu2014exact}), resulting in integer linear upper-level problems without coupling constraints and linear follower problems; see Appendix~\ref{app:detailed-instances} for details. We evaluate several behaviors. For strong-weak fixed measures, we have ${\beta \in \{0,0.3,0.5,0.7,1\}}$. For strong-weak endogenous measures, we test all $\beta(\cdot)\in\mathcal{B}$ with $\delta \in \{0.5,2,5,10\}$ for $\beta_{\text{Thr}}(\cdot)$, $a \in \{0.5,2,5,10\}$ for $\beta_{\text{Str}}(\cdot)$ and $\beta_{\text{Frg}}(\cdot)$, and $p \in \{2,5,10\}$ for $\beta_{\text{Str-p}}(\cdot)$ and $\beta_{\text{Frg-p}}(\cdot)$. For generalized measures, we consider all $\pi(\cdot)\in\Pi$, with $\gamma \in \{0.5,2,5,10,20,40,80\}$ for $\pi_{\text{Frg}}(\cdot)$ and $\gamma \in \{1,2,5,10,20\}$ for $\pi_{\text{Frg}^{+}}(\cdot)$ and $\pi_{\text{Str}}(\cdot)$. Finally, we set $\varepsilon = 0.01\cdot(\varphi_{ub}-\varphi_{lb})$.

As the proposed instances have a linear lower level, we reformulate the {\SAA} and {\DE} programs as single-level programs using the {\KKT} conditions, with complementarity constraints handled via SOS1 constraints. All optimization models were implemented in C++ and solved single-threaded with Gurobi 12.0.0 on an AMD EPYC 7532 CPU. To improve {\MCMC} Hit-and-Run mixing, the polyhedral set $\mathcal{S}(x,\varepsilon)$ was rounded using PolyRound 0.4.0~\citep{theorell2022polyround}. For the {\DE} {\SWDBO} baseline, nonlinear $\beta(\cdot)$ functions were approximated using PiecewiseLinApprox~\citep{codsi2025lina} with a $10^{-4}$ error tolerance. For generalized measures, the estimator $\hat{g}^{\mu}_{N'}(x)$ is computed in the second {\SAA} step using {\MCMC} Hit-and-Run samples of $\boldsymbol{\hat{y}}$ over $\mathcal{W}_{N'}$. All instances, results, and code are publicly available (\url{https://github.com/mariabazotte/endog-follower}).

\subsection{Computational Performance}

\subsubsection{Baseline Comparison for Strong-Weak Decision-Dependent Measures ($\mu_{\beta(\cdot)}$).}\label{subsubsec:baseline-performance}

\begin{figure}
    \FIGURE
    {
    \subcaptionbox{Time (seconds).\label{subfig:swd-saa-de-time}}
    {\includegraphics[width=0.38\textwidth]{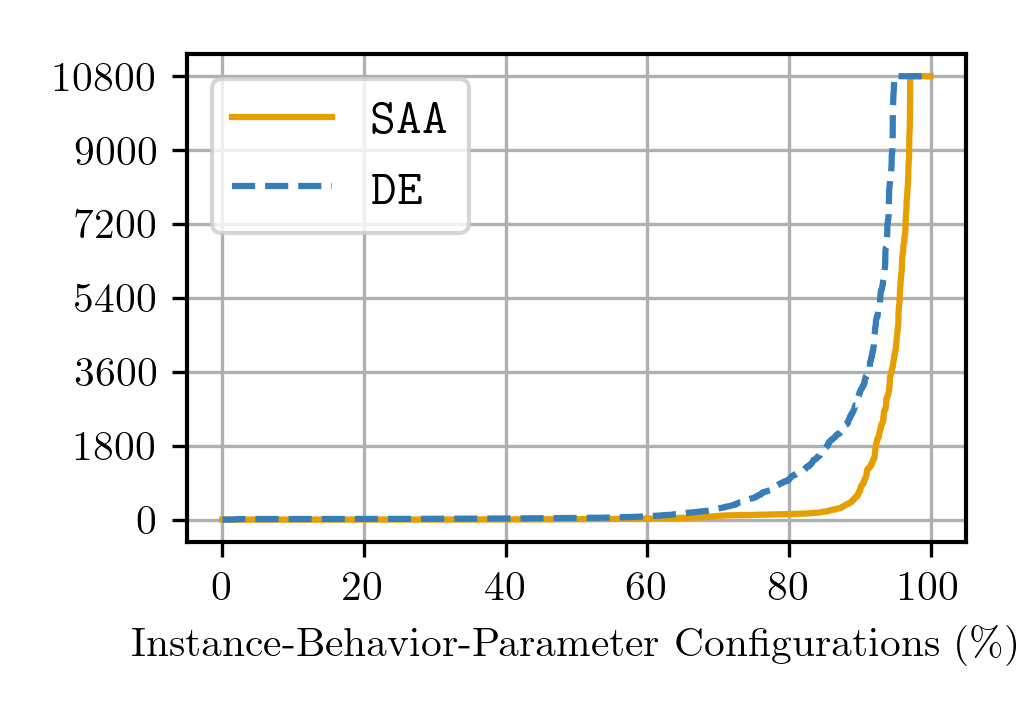}}
    \hfill
    \subcaptionbox{Optimality gap [$\text{gap}_{opt}]$($\%$). \label{subfig:swd-saa-de-gap}}
    {\includegraphics[width=0.38\textwidth]{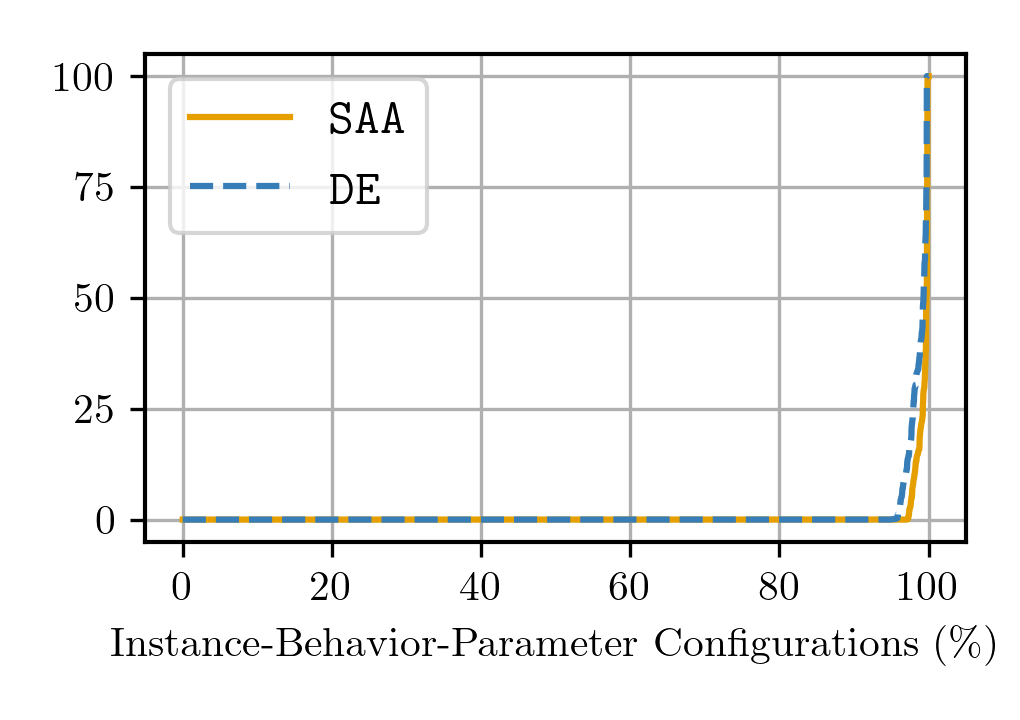}}
    }
    {Computational performance of {\SAA} Program~\eqref{Prog:saa-transformed-swd} and {\DE} Program~\eqref{Prog:de-swd} for $\mu_{\beta(\cdot)}$ measures.\label{fig:swd-saa-de}}
    {}
\end{figure}

We first benchmark our {\SAA} Program~\eqref{Prog:saa-transformed-swd} against the {\DE} Program~\eqref{Prog:de-swd}, as the latter is available in closed form for dichotomous measures. Figure~\ref{fig:swd-saa-de} reports the runtime and optimality gaps ($\text{gap}_{opt}$) for all configurations under a three-hour limit, using $N\!=\!100$ scenarios for the {\SAA} method. The latter is because the optimal value of the resulting {\SAA} program closely approximates the true objective value of the leader solution $\overline{x}$ returned by the method: the gap $[ 100\% \cdot (H(\overline{x}) + g^\mu(\overline{x}) - ub)/\lvert H(\overline{x}) + g^\mu(\overline{x}) \rvert]$ lies between $-1\%$ and $1\%$ for nearly $97\%$ of the configurations. Results are shown per percentage of runs, i.e., instance, endogenous follower behavior type $\beta(\cdot) \in \mathcal{B}$, and corresponding parameter-value ($\delta$, $a$, and $p$) configurations.

Figure~\ref{fig:swd-saa-de} shows that the proposed {\SAA} method outperforms the single-level {\DE} program in both optimality gaps and runtime across the considered configurations. Within the three-hour time limit, the single-level {\SAA} and {\DE} programs were solved to optimality for 97.1\% and 94.9\% of the instance-behavior-parameter configurations, respectively. For the remaining cases, the {\SAA} program produced optimality gaps in $(0,20\%]$, $(20\%,60\%]$, and above $60\%$ for 1.8\%, 0.9\%, and 0.2\% of the configurations, respectively, compared with 2.8\%, 1.8\%, and 0.5\% for the {\DE} program. The runtime results further highlight the superior performance of the proposed method. The {\SAA} program was solved in under 200 seconds for 85.4\% of the configurations, within 200--3600 seconds for 9.1\%, within 3600--7200 seconds for 2.0\%, and required more than 7200 seconds for 3.5\%. In comparison, the {\DE} program was solved within the same intervals for 67.6\%, 23.5\%, 2.8\%, and 6.1\% of the configurations, respectively. Overall, the {\SAA} program achieved an average runtime of 597.6 seconds versus 1081.7 seconds for the {\DE} program, corresponding to a 1.8x speedup, while also obtaining a lower average optimality gap (0.7\% versus 1.3\%).

We analyse the upper bound gap [\({{\text{gap}_{ub} = 100\% \cdot \left(v^{\mu}(\overline{x}^{\mu}) - v^{\mu}(\overline{x}^{\mu}_N)\right) / \lvert v^{\mu}(\overline{x}^{\mu} ) \rvert}}\)], where $\overline{x}^{\mu}$ and $\overline{x}^{\mu}_N$ denote the optimal or best feasible solutions obtained by the single-level {\DE} and {\SAA} programs, respectively. Table~\ref{tab:swd-saa-de-ubgap} reports this gap by parameterized behavior. The gap is consistently tight, averaging just $0.6\%$ across all configurations. Because the {\SAA} method yields leader solutions that have equal or superior quality to the solutions of the {\DE} in less time, we solely utilize it for endogenous measures in the subsequent analyses (see Appendix~\ref{app:additional-results-perf} for detailed bounds).

\begin{table}
\TABLE
{Average upper bound gap [$\text{gap}_{ub}$] ($\%$) per strong-weak decision-dependent type. \label{tab:swd-saa-de-ubgap}}
{
\begin{adjustbox}{width=0.62\textwidth,center}
\begin{tabular}{c|cc|cc|cc|cc|cc}
    \toprule
    $\beta_{\mathrm{Prp}}(\cdot)$ 
    & $\gamma$ & $\beta_{\mathrm{Thr}}^{\gamma}(\cdot)$ 
    & $a$ & $\beta_{\mathrm{Str}}^{a}(\cdot)$ 
    & $a$ & $\beta_{\mathrm{Frg}}^{a}(\cdot)$
    & $p$ & $\beta_{\mathrm{Str-p}}^{p}(\cdot)$
    & $p$ & $\beta_{\mathrm{Frg-p}}^{p}(\cdot)$ \\
    \midrule
    0.1$\pm$0.2 & 0.5 & -0.4$\pm$0.5 & 0.5 & 0.2$\pm$0.2 & 0.5 & 0.1$\pm$0.2 & 2 & 1.1$\pm$1.4 & 2 & 0.2$\pm$0.2 \\
    & 2 & -0.2$\pm$0.4 & 2 & 0.1$\pm$0.2 & 2 & 0.1$\pm$0.2 & 5 & 0.8$\pm$0.6  & 5 & 0.0$\pm$0.0 \\
    & 5 & 0.1$\pm$0.2 & 5 & 2.1$\pm$2.0 & 5 & 0.0$\pm$0.0 & 10 & 1.2$\pm$1.1 & 10 & 0.0$\pm$0.0 \\
    & 10 & 1.2$\pm$1.6 & 10  & 1.4$\pm$1.4  & 10 & 0.0$\pm$0.0 & & &\\
    \bottomrule
\end{tabular}
\end{adjustbox}
}
{\centering The values are averaged over all instances, with a 95\% confidence interval.}
\end{table}

\subsubsection{Comparison Across Decision-Dependent Measures.}\label{subsubsec:performance-all-measures}

We compare the performance of the {\SAA} method across the proposed measures in Section~\ref{sec:modeling}, and benchmark it against the commonly studied optimistic, pessimistic, and strong--weak fixed measures.

\begin{figure}
    \FIGURE
    {
    {\includegraphics[width=0.325\textwidth]{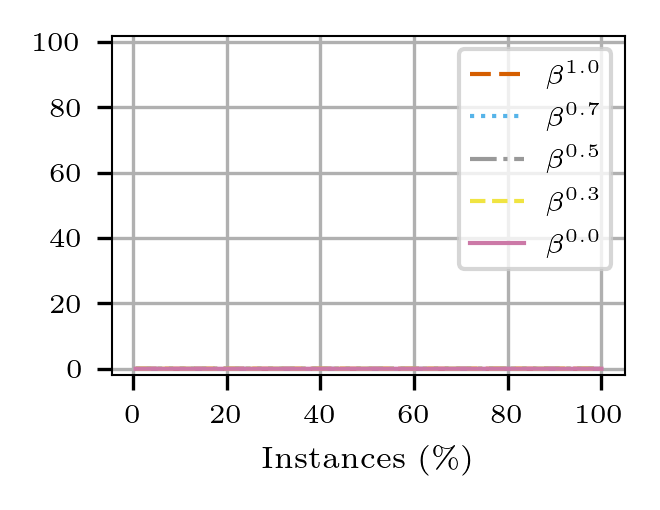}}
    \hfill
    {\includegraphics[width=0.325\textwidth]{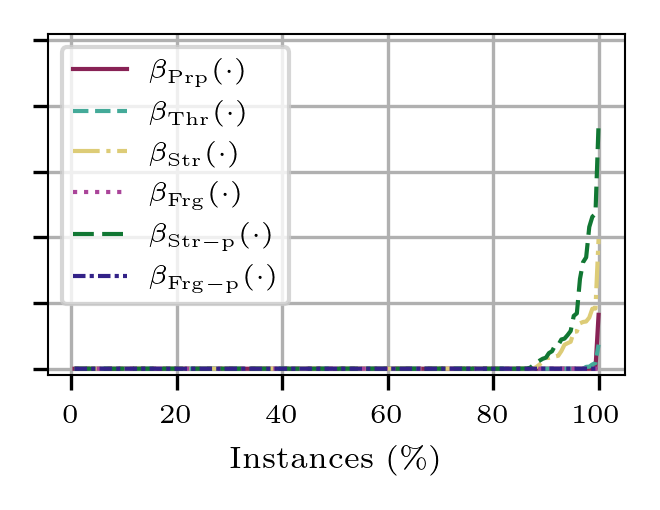}}
    \hfill
    {\includegraphics[width=0.325\textwidth]{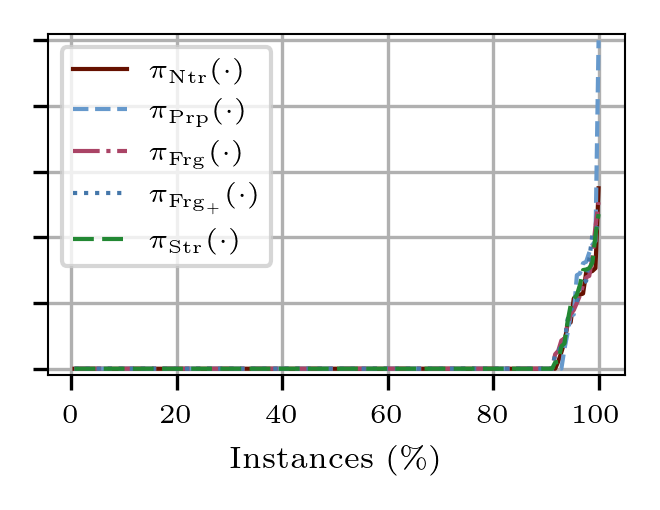}}
    }
    {Average optimality gap [$\text{gap}_{\text{opt}}$] (\%) by follower behavior type. \label{fig:comp-bv-gap}}
    {The values are averaged across parameterized variants of the same follower behavior (configurations), computed per instance.}
\end{figure}

\begin{figure}
    \FIGURE
    {
    {\includegraphics[width=0.325\textwidth]{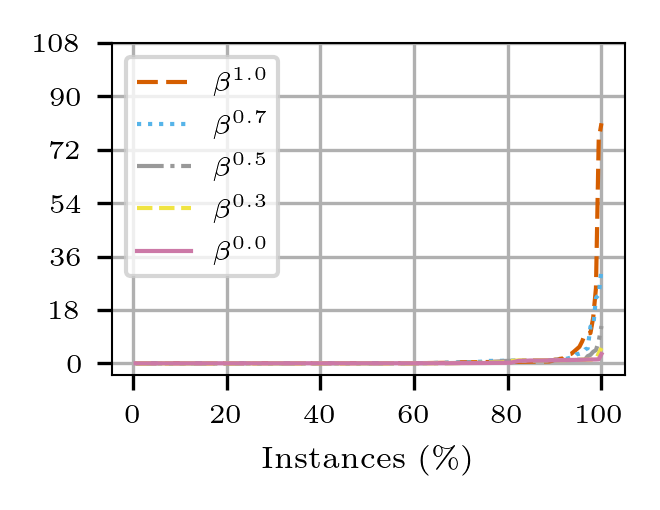}}
    \hfill
    {\includegraphics[width=0.325\textwidth]{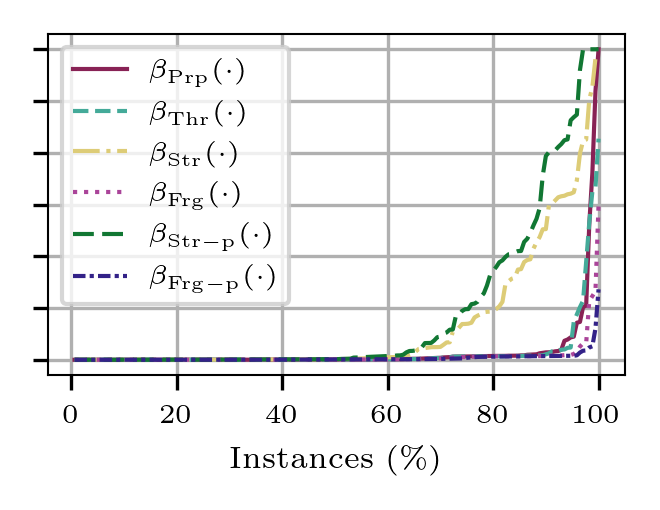}}
    \hfill
    {\includegraphics[width=0.325\textwidth]{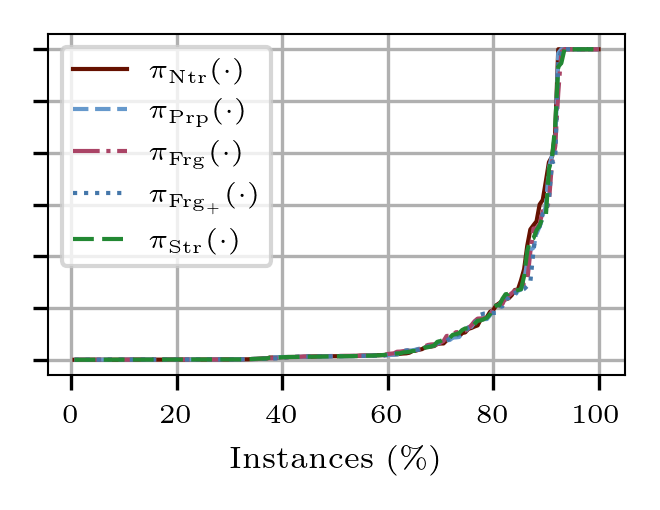}}
    }
    {Average computational time ($10^2$ seconds) by follower behavior type. \label{fig:comp-bv-time}}
    {The values are averaged across parameterized variants of the same follower behavior (configurations), computed per instance.}
\end{figure}

Figures~\ref{fig:comp-bv-gap} and~\ref{fig:comp-bv-time} present the results of all runs per type of follower behavior with a three-hour time limit. Figure~\ref{fig:comp-bv-time} reports the runtime in seconds, while Figure~\ref{fig:comp-bv-gap} shows the optimality gap, computed as described in the previous section with respect to the optimal value of the single-level of the {\DE}~Program~\eqref{Prog:de-swd} with fixed $\beta$ for strong--weak fixed measures, the single-level of the {\SAA}~Program~\eqref{Prog:saa-transformed-swd} with $N\!=\!100$ scenarios for strong-weak decision-dependent measures, and the single-level of the {\SAA}~Program~\eqref{Prog:saa-transf-int-cuts} with $N\!=\!1500$ scenarios (from a single Markov chain with thinning factor of ten) for generalized decision-dependent measures, per percentage of instances. For the {\SAA} Program~\eqref{Prog:saa-transf-int-cuts}, we first solve a separate optimization problem to determine $F_{lb}$ (see Section~\ref{subsec:methodology-ibl}) with a ten-minute time limit, taking the best lower bound if the problem is not solved to optimality; this preprocessing time is included in the total runtime of the {\SAA} program.

In these figures, $\beta^c$ denotes the measure $\mu_{\beta^c}$ with cooperation level $\beta \!=\! c$, where $\beta^1$ and $\beta^0$ correspond to the optimistic $\mu_{\Snew}$ and pessimistic $\mu_{\W}$ measures, respectively. Moreover, $\beta(\cdot) \in \mathcal{B}$ denotes the measure $\mu_{\beta(\cdot)}$, while $\pi(\cdot) \in \Pi$ denotes the measure $\mu_{\pi(\cdot)}$.

As expected, Figures~\ref{fig:comp-bv-gap} and~\ref{fig:comp-bv-time} show that fully endogenous models are computationally more challenging than support-endogenous ones. Within the time limit, the single-level {\DE} Program~\eqref{Prog:de-swd} for measures $\mu_\beta$ was solved to optimality for all instances and cooperation levels. The optimistic and pessimistic approaches exhibit similar performance, although a small subset of instances is solved faster under the pessimistic approach. This may be due to the use of $\varepsilon=0.01$, which allows near-optimal follower responses and can alter the degree of cooperation or adversarial behavior.

For the strong-weak decision-dependent measures $\mu_{\beta(\cdot)}$, the performance of the single-level {\SAA} Program~\eqref{Prog:saa-transformed-swd} depends strongly on the cooperation function $\beta(\cdot)$. Fragile behaviors, $\beta_{\text{Frg}}(\cdot)$ and $\beta_{\text{Frg-p}}(\cdot)$, yielded the best results: all configurations were solved to optimality, with over 95\% solved in under 200 seconds. Proportional and threshhold behaviors, $\beta_{\text{Prp}}(\cdot)$ and $\beta_{\text{Thr}}(\cdot)$, were also tractable, solving above 99\% of the configurations to optimality and nearly 90\% of configurations within 200 seconds. In contrast, sturdy behaviors, $\beta_{\text{Str}}(\cdot)$ and $\beta_{\text{Str-p}}(\cdot)$, were substantially more challenging, solving only 94.1\% and 90.8\% of the configurations to optimality, respectively, with 9\% and 14.9\% of configurations reaching the time limit. These results indicate that sturdy behaviors are significantly harder to solve than the other cooperation types. This indicates that the computational performance of the {\SAA} Program~\eqref{Prog:saa-transformed-swd} is influenced by the generated scenarios, i.e., the distribution of $\beta^{-1}(\zeta_\omega)$ for $\omega \in \mathcal{W}_N$, and thus by the structure of $\beta(\cdot)$. This is reinforced by the similar pattern observed for the {\DE} Program~\eqref{Prog:de-swd} (see Appendix~\ref{app:additional-results-perf}).

In contrast, the single-level {\SAA} Program~\eqref{Prog:saa-transf-int-cuts} exhibits similar performance across all generalized measures $\pi(\cdot)\in\Pi$. It was solved to optimality for 92.1\% of the configurations, while 1.1\%, 5.6\%, and 1.2\% had optimality gaps in $(0,10\%]$, $(10\%,35\%]$, and above $35\%$, respectively. Moreover, 60.8\% of the configurations were solved in under 200 seconds, 25.7\% within 200--3600 seconds, and 13.5\% required more than 3600 seconds, with 8.0\% reaching the time limit. In this case, computational performance is primarily driven by instance size: the number of binary leader variables affects the number of integer {\sf L}-shaped cuts, while the number of follower variables determines the dimension of the Hit-and-Run Markov chain used during cut separation (see Appendix~\ref{app:additional-results-perf}).

The single-level {\SAA} Program~\eqref{Prog:saa-transf-int-cuts} for measures $\mu_{\pi(\cdot)}$ required longer runtimes than the {\SAA} Program~\eqref{Prog:saa-transformed-swd} for measures $\mu_{\beta(\cdot)}$. This is primarily due to the use of {\MCMC} Hit-and-Run sampling during cut separation, as the expectation under generalized measures has no closed-form expression. Moreover, {\SAA} Program~\eqref{Prog:saa-transf-int-cuts} uses a binary representation of the leader variables to derive integer {\sf L}-shaped cuts, resulting in larger problem sizes. Furthermore, the evaluation of the resulting solution $\overline{x}_N^\mu$ also differs between the two classes of measures. For measures $\mu_{\beta(\cdot)}$, $g^\mu(\overline{x}_N^\mu)$ can be computed exactly from the enumerable scenarios. For measures $\mu_{\pi(\cdot)}$, $\hat g^\mu_{N'}(\overline{x}_N^\mu)$ is estimated via Hit-and-Run using $N'=4\times10^4$ samples from four independent Markov chains with thinning factor 50. The resulting chains exhibited potential scale reduction factors $\hat R$ close to one, providing evidence of convergence~\citep{gelman1992inference}. The evaluation procedure completed in under 60 seconds for 74.3\% of the configurations and exceeded 1800 seconds for only 5.0\%.

\subsection{Comparison of Follower's Behavior}\label{subsec:comparison-behaviors}

\subsubsection{Leader and Follower Objective Values.}\label{subsubsec:objective-values}

We compare the leader and follower expected objective values ($\mathbb{E}_{\boldsymbol{\hat{y}}\!\sim\!\mu_x}\![F(x,\boldsymbol{\hat{y}})]$ and $\mathbb{E}_{\boldsymbol{\hat{y}}\!\sim\!\mu_x}\![f(x,\boldsymbol{\hat{y}})]$) of the optimal or best leader solutions obtained with the methods in Section~\ref{subsubsec:performance-all-measures} for distinct follower behaviors within the three-hour time limit. We use the leader and follower objective values under the optimistic approach as references, and focus our analysis on behaviors that exhibit the most distinct patterns. To avoid biases, we restrict this analysis to instances where the {\DE} or {\SAA} program achieved an optimality gap below 5\%.

For a given behavior or measure $\mu$, let $v^\mu$ and $v_f^\mu$ denote the expected leader and follower objective values of the optimal or best leader solution under $\mu$. For measures ${\mu \! = \! \mu_\beta}$, we set ${v^\mu \! = H\left(\overline{x}^\mu\right) + g^\mu\!\left(\overline{x}^\mu\right)}$ and ${v_f^\mu \! = g_f^\mu \left( \overline{x}^\mu \right)}$ with $\overline{x}^\mu$ obtained from the single-level {\DE}~Program~\eqref{Prog:de-swd}; for measures $\mu \! = \! \mu_{\beta(\cdot)}$, we set ${v^\mu \! = H\!\left(\overline{x}_N^\mu\right) + g^\mu\!\left(\overline{x}_N^\mu\right)}$ and ${v_f^\mu \! = g_f^\mu\!\left(\overline{x}_N^\mu\right)}$ with $\overline{x}_N^\mu$ obtained from the single-level {\SAA}~Program~\eqref{Prog:saa-transformed-swd}; and for measures $\mu \! = \! \mu_{\pi(\cdot)}$, we set ${v^\mu \! = H\!\left(\overline{x}_N^\mu\right) + \hat{g}_{N'}^\mu\!\left(\overline{x}_N^\mu\right)}$ and ${v^\mu_f \! = \hat{g}_{f,N'}^\mu\!\left(\overline{x}_N^\mu\right)}$, with $\overline{x}_N^\mu$ obtained from the single-level {\SAA}~Program~\eqref{Prog:saa-transf-int-cuts}. The parameters $N$ and $N'$ are defined as in the previous section. The leader and follower expected reference objective gaps are ${\text{gap}^\mu \! = 100\% \cdot (v^{\beta^{1}} - v^\mu)/\lvert v^{\beta^{1}} \rvert}$ and ${\text{gap}_f^\mu \! = 100\% \cdot ( v_f^{\beta^{1}} - v_f^\mu)/\lvert  v_f^{\beta^{1}} \rvert}$, where $\beta^{1}$ is ${\mu_{\beta^{1}} = \mu_{\Snew}}$. 

\begin{figure}
    \FIGURE
    {
    {\includegraphics[height=4cm]{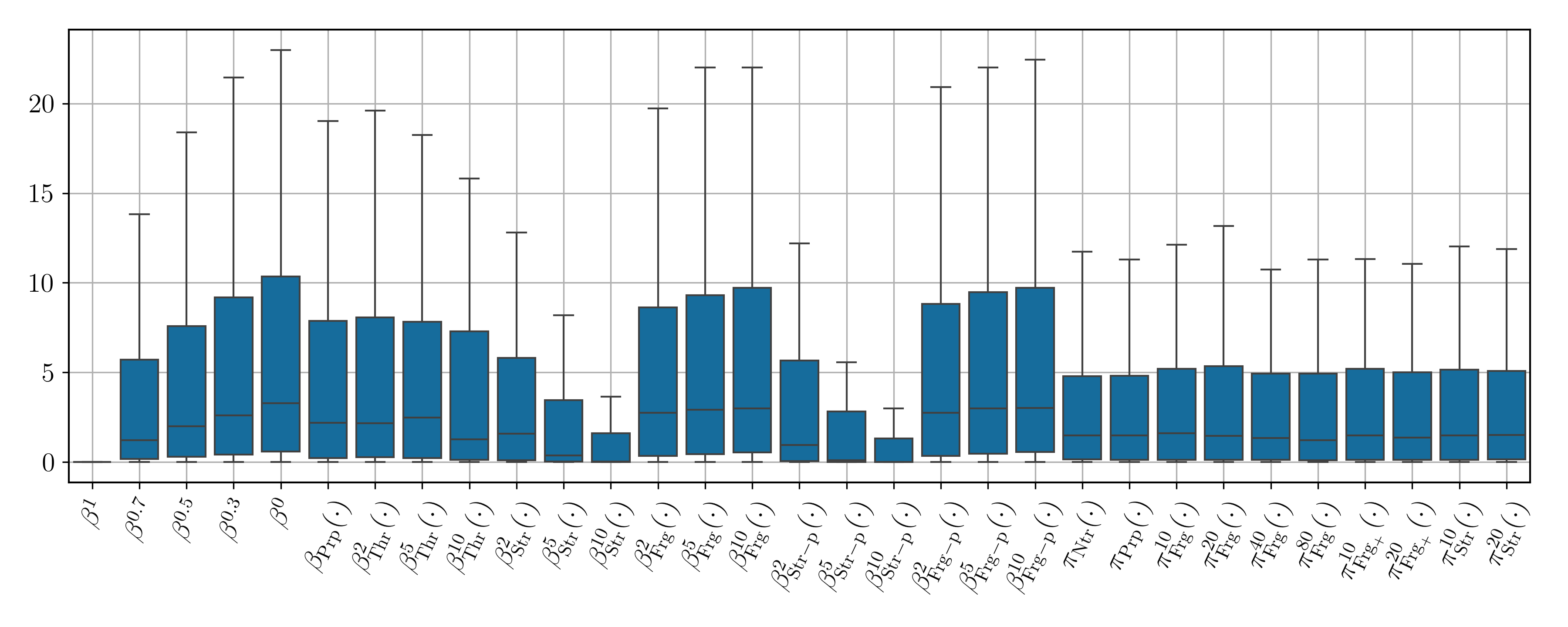}}
    }
    {Median and quantiles of the leader expected reference objective gap [$\text{gap}^\mu$] ($\%$). \label{fig:leader-reference-gap}}
    {}
\end{figure}

\begin{figure}
    \FIGURE
    {
    {\includegraphics[height=4cm]{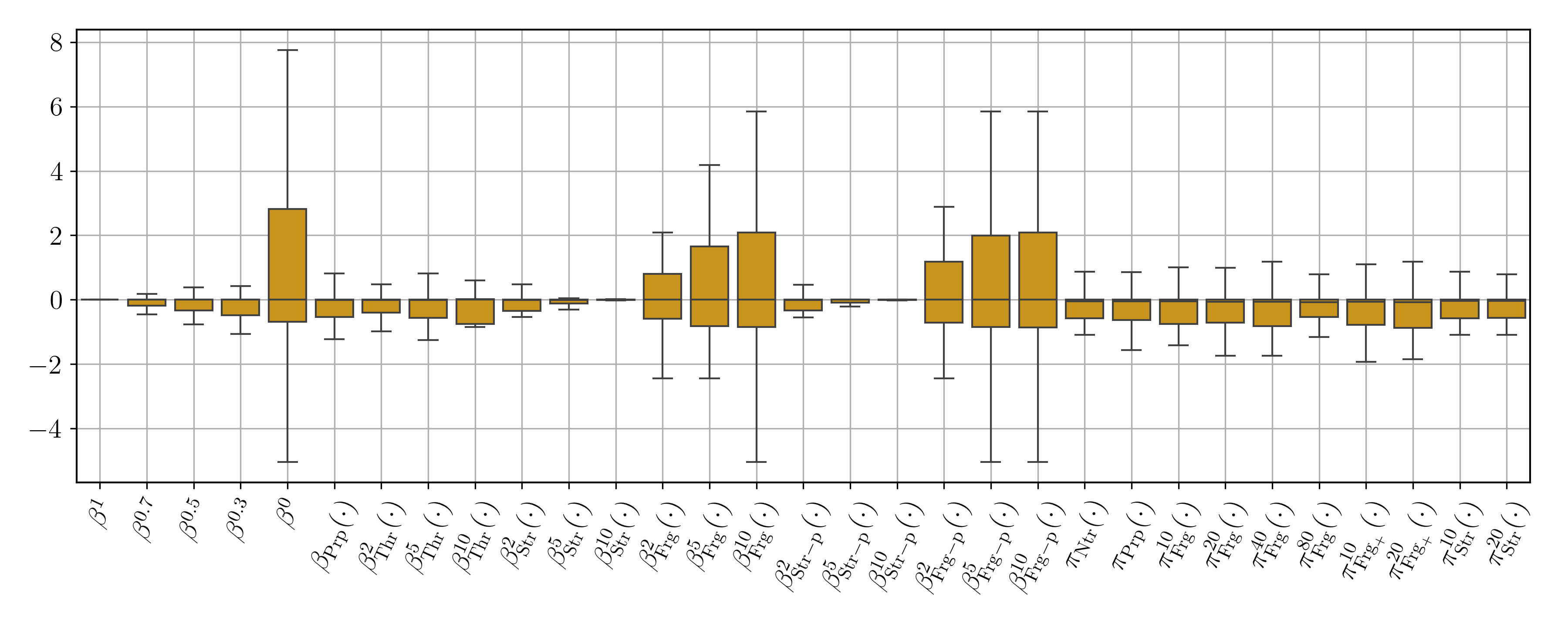}}
    }
    {Median and quantiles of the follower expected reference objective gap [$\text{gap}_f^\mu$] ($\%$) [All instances]. \label{fig:follower-reference-gap}}
    {}
\end{figure}

Figures~\ref{fig:leader-reference-gap} and~\ref{fig:follower-reference-gap} report the median and quantiles of the leader and follower expected reference objective gap for each behavior. We use the notation $\beta_{\text{Thr}}^\delta(\cdot)$, $\beta_{\text{Str}}^a(\cdot)$, $\beta_{\text{Frg}}^a(\cdot)$, $\beta_{\text{Str-p}}^p(\cdot)$, $\beta_{\text{Frg-p}}^p(\cdot)$, $\pi_{\text{Frg}}^\gamma(\cdot)$, $\pi_{\text{Frg}^{+}}^\gamma(\cdot)$, and $\pi_{\text{Str}}^\gamma(\cdot)$, where $\delta$, $a$, $p$, and $\gamma$ denote the parameter values of the behaviors.

Figure~\ref{fig:leader-reference-gap} illustrates the impact of follower behavior on the leader’s expected objective value for the optimal solutions of the {\DE} or {\SAA} programs. As discussed in Section~\ref{subsec:literature-background}, the optimistic and pessimistic objective values provide lower and upper bounds, respectively, for more general follower behaviors. The pessimistic behavior yielded the highest average leader expected reference gap, namely 9.0\%. Moreover, the results support Proposition~\ref{prop:swd}. Sturdy behaviors $\beta_{\text{Str}}$ and $\beta_{\text{Str-p}}$ produced average gaps ranging from 2.5\% to 5.3\% and 2.4\% to 5.2\%, respectively, while the proportional behavior yielded an average gap of 6.5\%. In contrast, fragile behaviors $\beta_{\text{Frg}}$ and $\beta_{\text{Frg-p}}$ exhibited larger gaps, ranging from 7.2\% to 7.9\%.

For the generalized endogenous behaviors, the average leader reference gap remains stable across configurations, at around 4.9--5.3\%. This suggests that, although the pessimistic and optimistic solutions can differ more significantly in leader value, the set of follower $\varepsilon$-optimal solutions $\mathcal{S}(x,\varepsilon)$ remains relatively small for the considered instances. Hence, the gaps obtained under general endogenous measures over the entire set $\mathcal{S}(x,\varepsilon)$ are similar to those of the neutral approach.

Conversely, Figure~\ref{fig:follower-reference-gap} shows that, for most behaviors, the follower’s expected reference gap is centered around zero and takes negative values, indicating that the follower’s expected objective value is preserved or improved under these behaviors, while the leader’s objective deteriorates (see Figure~\ref{fig:leader-reference-gap}). Highly adversarial behaviors, such as the pessimistic and strong--weak fragile behaviors, however, can produce positive reference gaps, worsening the follower's expected outcome. Similar trends are observed when considering the follower’s optimal objective $\varphi(x)$ instead of the expected objective $\mathbb{E}_{\boldsymbol{\hat{y}}\sim\mu_x}[f(x,\boldsymbol{\hat{y}})]$; see Appendix~\ref{app:additional-results-comp}.

\begin{figure}
    \FIGURE
    {
    {\includegraphics[height=4cm]{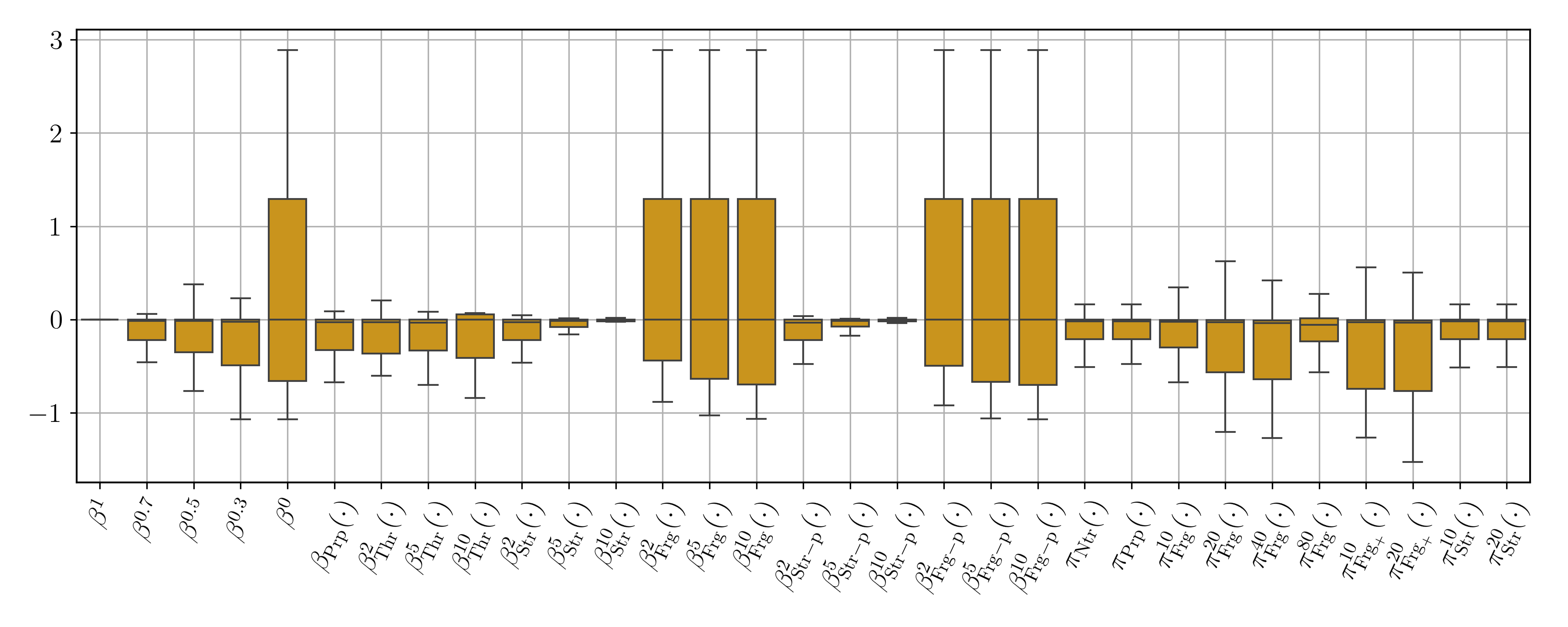}}
    }
    {Median and quantiles of the follower expected reference objective gap [$\text{gap}_f^\mu$] ($\%$) [$\alpha \in [-1,-0.1)$]. \label{fig:follower-reference-gap-alpha-negative}}
    {}
\end{figure}

\begin{figure}
    \FIGURE
    {
    {\includegraphics[height=4cm]{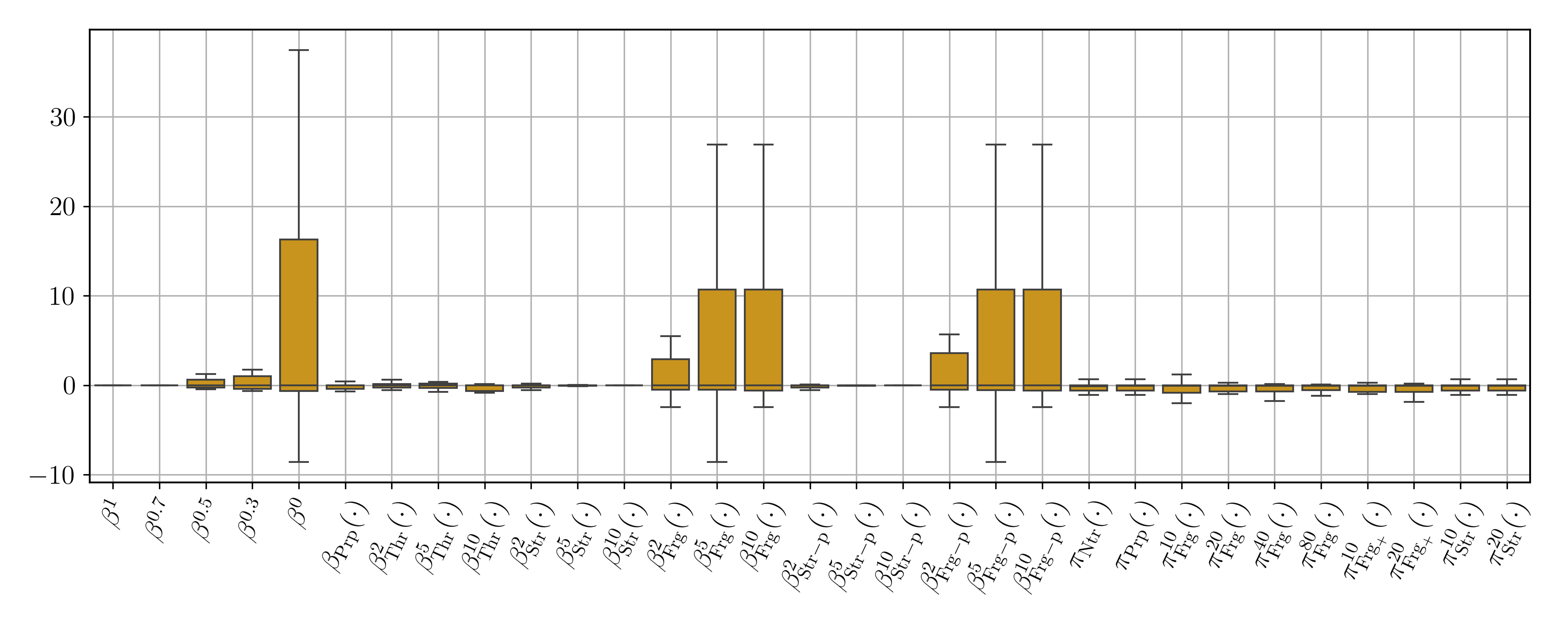}}
    }
    {Median and quantiles of the follower expected reference objective gap [$\text{gap}_f^\mu$] ($\%$) [$\alpha \in (-0.1,0.1)$]. \label{fig:follower-reference-gap-alpha-middle}}
    {}
\end{figure}

\begin{figure}
    \FIGURE
    {
    {\includegraphics[height=4cm]{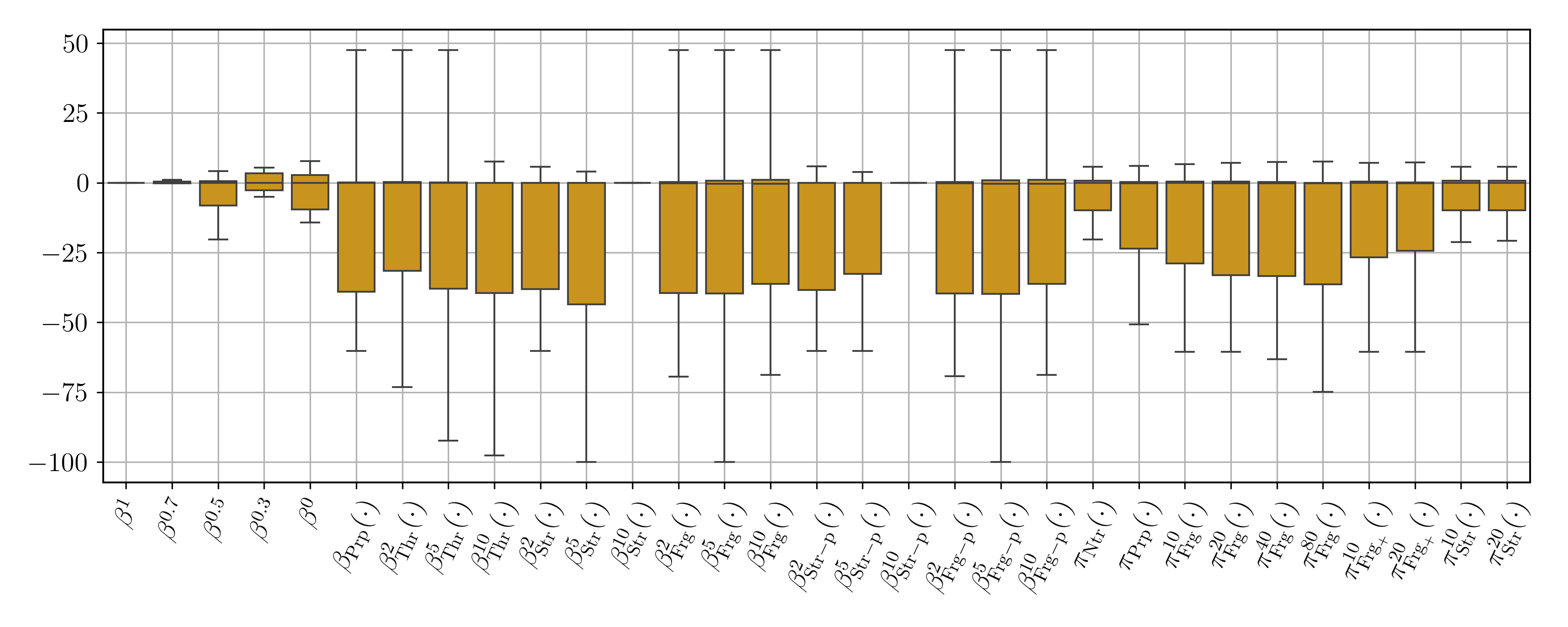}}
    }
    {Median and quantiles of the follower expected reference objective gap [$\text{gap}_f^\mu$] ($\%$) [$\alpha \in (0.1,1]$]. \label{fig:follower-reference-gap-alpha-positive}}
    {}
\end{figure}

To analyze this further, we study the follower’s gap as a function of the objective alignment $\alpha = d^\top e / (\|d\|_2 \|e\|_2) \in [-1,1]$ with respect to $\hat{y}$, where $f_{\text{Lin}}(x,\hat{y}) = e^\top \hat{y}$ and $F_{\text{Lin}}(x,\hat{y}) = d^\top \hat{y}$ are the linear objective components. We partition the instances into three ranges: fully opposed to weakly opposed ($\alpha \in [-1,-0.1)$, 71 instances, Figure~\ref{fig:follower-reference-gap-alpha-negative}); weakly aligned/orthogonal ($\alpha \in [-0.1,0.1)$, 56 instances, Figure~\ref{fig:follower-reference-gap-alpha-middle}); and aligned ($\alpha \in [0.1,1]$, 43 instances, Figure~\ref{fig:follower-reference-gap-alpha-positive}).


Across all alignment ranges and behaviors, the median gap remains close to zero. For instances with more aligned objectives, $\alpha \in (0.1,1]$, the follower’s expected reference objective gap is predominantly negative and can attain substantially negative values across all behaviors, indicating that more general or manipulative behaviors may significantly improve the follower’s expected objective relative to the optimistic behavior. In contrast, for more opposed objectives, $\alpha \in [-1,-0.1)$, the gaps remain much closer to $0\%$, typically between $-2\%$ and $3\%$, although most behaviors still exhibit predominantly negative gaps. Highly adversarial behaviors, however, produce both positive and negative gaps, suggesting a less predictable impact that depends more strongly on instance characteristics. For nearly independent objectives, $\alpha \in [-0.1,0.1)$, most behaviors again yield gaps close to zero, while highly adversarial behaviors can generate substantially positive gaps, indicating that weakly correlated objectives may amplify the variability and detrimental effects of adversarial behavior on the follower’s expected objective value.

\subsubsection{Leader Objective Value Loss under Misspecified Follower Behavior.}\label{subsubsec:misspecification}

We evaluate the loss in the leader’s objective due to follower behavior misspecification. Let $v^{\mu \leftarrow \tilde{\mu}}$ denote the leader's expected objective value evaluated under the follower behavior or endogenous measure $\mu$, using the optimal leader decision obtained under a distinct follower behavior or measure $\tilde{\mu} \neq \mu$. For measures $\mu \! = \! \mu_\beta$ and $\mu \! = \! \mu_{\beta(\cdot)}$, we define $v^{\mu \leftarrow \tilde{\mu}} \! = H\!\left(x^{\tilde{\mu}}\right) + g^\mu\!\left(x^{\tilde{\mu}}\right)$, and for measures $\mu = \mu_{\pi(\cdot)}$, we define $v^{\mu \leftarrow \tilde{\mu}} \! = H\!\left(x^{\tilde{\mu}}\right) + \hat{g}_{N'}^\mu\!\left(x^{\tilde{\mu}}\right)$. Here, $x^{\tilde{\mu}}$ denotes the optimal or best leader decision as defined in Section~\ref{subsubsec:objective-values}. As in the previous section, we restrict our analysis to instances for which the single-level {\DE} or {\SAA} program achieved optimality gaps of at most 5\%. The leader ($\mu\leftarrow\tilde{\mu}$) expected misspecified objective gap $\text{gap}^{\mu\leftarrow\tilde{\mu}}=100\%\cdot(v^\mu - v^{\mu \leftarrow \tilde{\mu}})/\lvert v^\mu \rvert$, with $v^\mu$ as defined in Section~\ref{subsubsec:objective-values}. 

\begin{figure}
    \FIGURE
    {
    {\includegraphics[height=4cm]{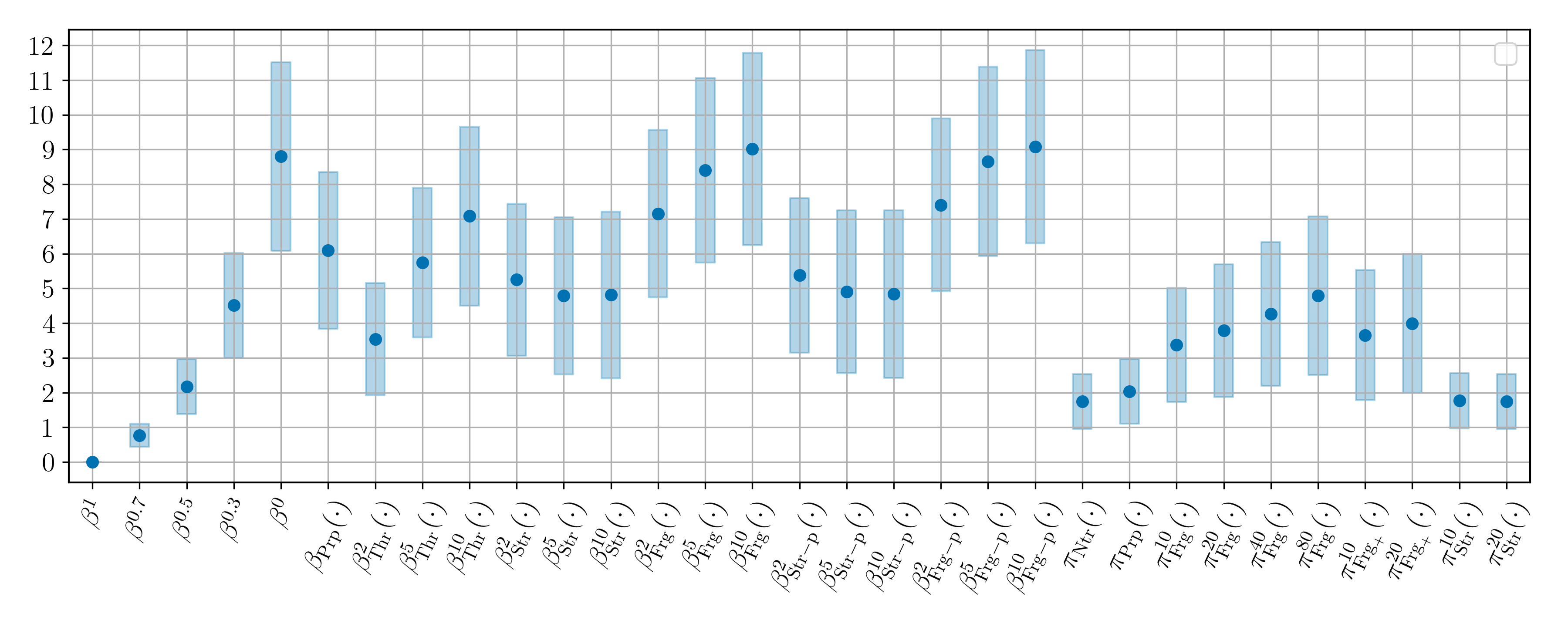}}
    }
    {Average leader ($\mu \leftarrow \beta^1$) expected misspecified objective gap [$\text{gap}^{\mu\leftarrow \beta^1}$] ($\%$). \label{fig:leader-misspecified-gap-optimistic}}
    {The values are averaged over all instances, and the shaded area represents the 95\% confidence interval.}
\end{figure}

\begin{figure}
    \FIGURE
    {
    {\includegraphics[height=4cm]{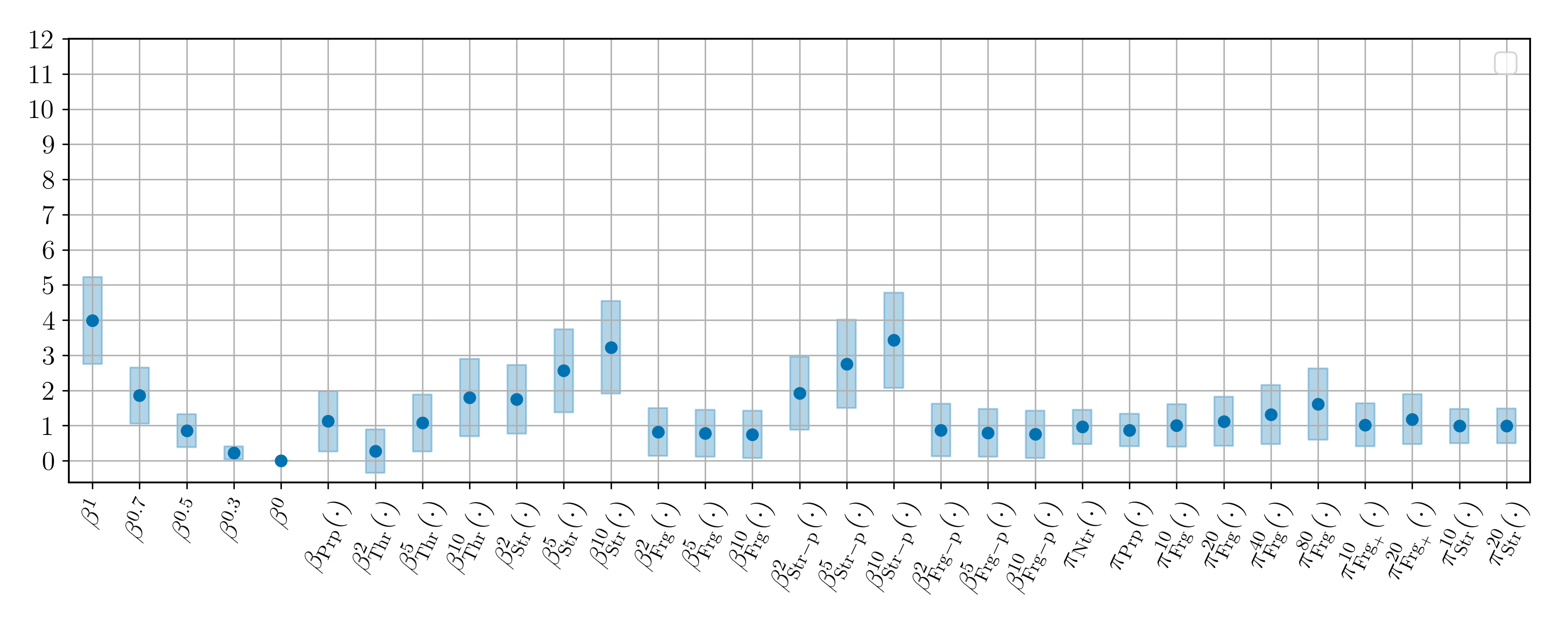}}
    }
    {Average leader ($\mu \leftarrow \beta^0$) expected misspecified objective gap [$\text{gap}^{\mu\leftarrow \beta^0}$] ($\%$). \label{fig:leader-misspecified-gap-pessimistic}}
    {The values are averaged over all instances, and the shaded area represents the 95\% confidence interval.}
\end{figure}

Figures~\ref{fig:leader-misspecified-gap-optimistic} and \ref{fig:leader-misspecified-gap-pessimistic} report the average leader $\text{gap}^{\mu \leftarrow \beta^1}$ and $\text{gap}^{\mu \leftarrow \beta^0}$ across instances.We focus on these approaches because they yield simpler and more tractable formulations and are the most widely studied in the literature. This allows us to assess whether considering only these canonical approaches is sufficient in practice. These figures show that misspecifying the follower’s behavior can lead to substantial losses in the leader’s objective value. Leader solutions based on standard follower behaviors may be suboptimal when evaluated against alternative behaviors, highlighting the importance for the leader of accurately modeling the follower behavior (e.g., by investigating historical data on the follower's responses) when making decisions.

\begin{figure}
    \FIGURE
    {
    {\includegraphics[height=4cm]{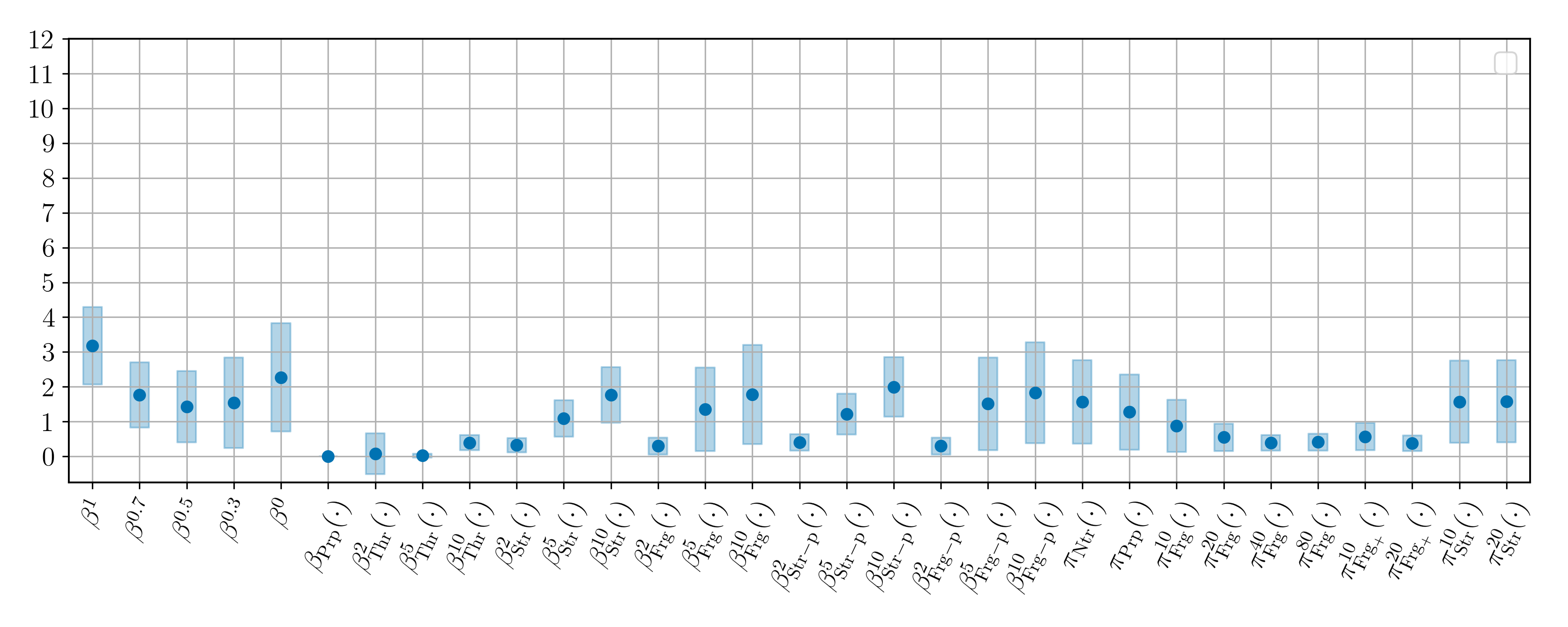}}
    }
    {Average leader ($\mu \leftarrow \beta_{\text{Prp}}(\cdot)$) expected misspecified objective gap [$\text{gap}^{\mu\leftarrow \beta_{\text{Prp}}(\cdot)}$] ($\%$). \label{fig:leader-misspecified-gap-proportional}}
    {The values are averaged over all instances, and the shaded area represents the 95\% confidence interval.}
\end{figure}

The $(\mu\leftarrow\beta^1)$ gaps for leader solutions obtained under the optimistic approach (Figure~\ref{fig:leader-misspecified-gap-optimistic}) are small for more cooperative follower behaviors, such as the strong-weak sturdy behaviors, and large for more adversarial behaviors, such as pessimistic or strong-weak fragile behaviors. The opposite pattern occurs for the $(\mu\leftarrow\beta^0)$ gap for the pessimistic approach (Figure~\ref{fig:leader-misspecified-gap-pessimistic}). The figures show that using leader decisions based on the optimistic approach leads to larger losses in the leader’s objective value when the follower adopts an alternative behavior, compared to decisions obtained from the pessimistic approach. Moreover, behaviors that induce a balanced distribution between pessimistic and optimistic outcomes tend to yield low leader misspecification gaps, as the sensitivity of the leader’s objective to the follower’s actual behavior is reduced. Under the strong–-weak proportional behavior, the average misspecified leader objective gap remains small for both more cooperative or adversarial types of behaviors (see Figure~\ref{fig:leader-misspecified-gap-proportional}). Hence, such balanced behaviors lead to more robust leader decisions, even though they may still incur substantial losses.

\paragraph{\textbf{Guidelines}} From the presented experimental results, we derive the following guidelines: (i) From the leader's perspective, it is crucial to correctly model follower behavior when making decisions, for instance by analyzing historical data on follower responses and, more generally, by understanding the follower’s tendencies. When such information is not available, it is advisable to adopt behaviors that balance pessimistic and optimistic outcomes, such as the strong–weak proportional $\beta_{\text{Prp}}(\cdot)$; (ii) From the follower's perspective, particularly when the alignment $\alpha$ between the leader’s and follower’s objectives lies in $(0.1, 1.0]$ (i.e., when the objectives are aligned), it is beneficial to consider diverse behaviors, as these can improve in the follower’s outcome.

\section{Conclusions}\label{sec:conclusions}

We investigated {\IBO}s in which the follower’s selected optimal response is uncertain and decision-dependent. We introduced endogenous measures capturing different behaviors, including strong--weak measures restricted to optimistic and pessimistic responses, and generalized measures supported on the entire optimal response set. We proposed the {\TIBO} reformulation that expresses the endogenous uncertainty as a transformation function of the leader’s decisions and additional exogenous randomness. We then presented an inverse transformation for dichotomous measures and a Markov chain-based transformation for generalized measures, and we developed {\SAA}-based solution methods. For dichotomous measures, we derived an {\SAA} program that avoids the nonlinearities induced by the endogenous probability measure in the {\DE} program. For generalized measures, where the endogenous expectation lacks a closed form, we combined integer {\sf L}-shaped cuts with {\MCMC} simulation to implicitly handle the transformation when leader decisions are binary.

Our experiments showed that the proposed {\SAA} program for dichotomous measures outperforms the {\DE} program. Despite the additional complexity introduced by simultaneously modeling endogenous probability and support, the proposed methodology remains computationally tractable. We also demonstrated that follower behavior misspecification can substantially affect both leader and follower outcomes and lead to suboptimal leader decisions. Several research directions remain open, including developing methods for generalized measures with more general leader decision spaces, extending the proposed transformations to integer follower responses, deriving performance guarantees under different follower behaviors, and investigating combinatorial follower structures.






\bibliographystyle{informs2014} 
\bibliography{refs}

\begin{APPENDICES}

\section{Discrete and Absolute Continuous Decision-Dependent Measures} \label{app:intermediate-discrete-continuous}

In this section, we develop the {\IBO} Program~\eqref{Prog:intermediate} for discrete and absolute continuous supports. When ${\mathcal{S}(x,\varepsilon)=\{\hat{y}_1(x,\varepsilon), \ldots, \hat{y}_{s(x,\varepsilon)}(x,\varepsilon)\}}$ for all $x \in \mathcal{X}$, i.e., the follower’s $\varepsilon$-optimal response set is finite and discrete, the leader specifies the likelihood of each $\varepsilon$-optimal follower response, assigning a discrete probability distribution over $\mathcal{S}(x,\varepsilon)$ for each $x \in \mathcal{X}$. The Discrete {\IBO} is:
\begin{equation}
    \Theta_{\I}(\mu_P) = \min_{x\in \mathcal{X}} \left\{ H(x) + \sum_{j = 1}^{s(x,\varepsilon)} p_{j}(x,\varepsilon) F(x,\hat{y}_j(x,\varepsilon)) \right\},
\end{equation}
where $p_j(x,\varepsilon)$ is the probability that the follower selects the $\varepsilon$-optimal solution ${\hat{y}_j(x) \in \mathcal{S}(x,\varepsilon)}$, and ${\mu_P=\cup_{x \in \mathcal{X}} \{p_1(x,\varepsilon),\ldots,p_{s(x,\varepsilon)}(x,\varepsilon)\}}$ represents the set of probability distributions for all leader decisions. For probability measures $\mu_x^{\pi(\cdot)}$ in Equation~\eqref{eq:measure-pi}, we have ${p_j(x,\varepsilon) =  \pi\!\left( \tilde{F}(x,\hat{y}_j), \tilde{\varphi}(x)\right)/\sum_{i = 1}^{s(x,\varepsilon)} \pi\!\left( \tilde{F}(x,\hat{y}_i), \tilde{\varphi}(x)\right)}$. 
Furthermore, when $\mu_x$ is absolutely continuous with respect to the Lebesgue measure for each $x \in  \mathcal{X}$, with density $\rho_x$, the Continuous {\IBO} is then given by:
\begin{equation}
    \Theta_{\I}(\mu_{\rho}) = \min_{x \in \mathcal{X}} \left\{ H(x) + \int_{\hat{y} \in \mathcal{S}(x,\varepsilon)} F(x,\hat{y}) \, \rho_x(\hat{y}) \, d\hat{y} \right\},
\end{equation}
where $\mathcal{S}(x,\varepsilon)$ is assumed to be Lebesgue measurable with positive measure for all $x \in \mathcal{X}$, and $\mu_\rho = \cup_{x \in \mathcal{X}} \{\mu_{\rho_x}\}$ denotes the collection of densities corresponding to all leader decisions. For $\mu_x^{\pi(\cdot)}$ in Equation~\eqref{eq:measure-pi}, we have ${\rho_x(\hat{y}) = \pi\left(\tilde{F}(x,\hat{y}), \tilde{\varphi}(x)\right)/
\int_{\mathcal{S}(x,\varepsilon)} \pi\left(\tilde{F}(x,z), \tilde{\varphi}(x)\right)\, d\hat{y}}$ for $\hat{y} \in \mathcal{S}(x,\varepsilon)$.

\section{Examples} \label{app:examples}

\begin{example}\label{ex:behavior-swd}
    Consider the {\SBO}, {\WBO}, and {\SWDBO} approaches for the following problem:
    \begin{equation}
        \mbox{{\SBO}:} \min_{0 \leq x \leq 2} \left\{ 19x + \min_{\hat{y} \in \mathcal{S}(x,0)} \left(10\hat{y}_1 + 20\hat{y}_2 \right)  \right\}, \nonumber
    \end{equation}
    \begin{equation}
        \mbox{{\WBO}:} \min_{0 \leq x \leq 2} \left\{ 19x + \max_{\hat{y} \in \mathcal{S}(x,0)} \left(10\hat{y}_1 + 20\hat{y}_2 \right)  \right\}, \nonumber
    \end{equation}
    \begin{equation}
        \mbox{{\SWDBO}:} \min_{0 \leq x \leq 2} \left\{ 19x + \beta \left( \tilde{\varphi}(x) \right) \min_{\hat{y} \in \mathcal{S}(x,0)} \left(10\hat{y}_1 + 20\hat{y}_2 \right)  + [ 1- \beta \left( \tilde{\varphi}(x) \right) ] \max_{\hat{y} \in \mathcal{S}(x,0)} \left(10\hat{y}_1 + 20\hat{y}_2 \right)  \right\}, \nonumber
    \end{equation}
    \begin{equation}
        \mbox{with }\mathcal{S}(x,0) = \arg\min_y \left\{ y_1 + y_2 : y_1 + y_2 \geq 3 - x, y_1 \geq 0, y_2 \geq 0 \right\}. 
    \end{equation}
    The follower’s optimal set reduces to ${\mathcal{S}(x,0) = \{y_1 \geq 0,\, y_2 \geq 0 : y_1 + y_2 = 3 - x\}}$, yielding the follower’s value function $\varphi(x) = 3 - x$ with bounds $\varphi_{lb} = 1$ and $\varphi_{ub} = 3$. Given a leader’s feasible decision $x$, the optimistic and pessimistic follower responses are ${\hat{y}^{\Snew}(x) = \{ \hat{y}_1^{\Snew}(x) = 3-x,\, \hat{y}_2^{\Snew}(x) = 0 \}}$ with $F_{{\Snew}}(x,\varepsilon) = 10(3-x)$ and ${\hat{y}^{\W}(x) = \{ \hat{y}_1^{\W}(x) = 0,\, \hat{y}_2^{\W}(x) = 3-x \}}$ with $F_{{\W}}(x,\varepsilon)=20(3-x)$, respectively. Thus, the {\SWDBO} reduces to:
    \begin{equation}
        \mbox{{\SWDBO}:} \min_{0 \leq x \leq 2} \left\{ 60 - x - 10 (3-x) \beta \left( \tilde{\varphi}(x) = 1 - \tfrac{x}{2}  \right) \right\},  \nonumber
    \end{equation}
    and the cooperation level functions are defined as:
    \begin{subequations}
        \begin{align*}
            \beta_{\text{Prp}}\left( \tilde{\varphi}(x) \right) &= \tfrac{x}{2},& \\
            \beta_{\text{Thr}}\left( \tilde{\varphi}(x) \right) &= \tfrac{1}{1+\exp \left( \frac{ \delta (1-x) }{2} \right)},& \delta > 0, \\
            \beta_{\text{Str}}\left(\tilde{\varphi}(x)\right) &= \tfrac{1- \exp\left(\frac{-a x}{2} \right)}{1- \exp(-a)},& a>0, \\
            \beta_{\text{Frg}}\left(\tilde{\varphi}(x) \right) &= 1 - \tfrac{1-\exp \left(-a\left[1-\frac{x}{2}\right] \right)}{1-\exp(-a)},& a>0, \\
            \beta_{\text{Str-p}}\left(\tilde{\varphi}(x)\right) &= 1 - \left(1 -\tfrac{x}{2} \right)^p,& p>1, \\
            \beta_{\text{Frg-p}}\left(\tilde{\varphi}(x) \right) &= \left(\tfrac{x}{2} \right)^p,& p>1.
        \end{align*}
    \end{subequations}
    Table~\ref{tab:example-swd-results} reports the leader’s optimal solution $x_{\mu}^*$ and optimal objective $\Theta_{\I}(\mu)$ for each measure $\mu$, as well as the follower's optimal value given the leader's solution $\varphi(x^{*}_{\mu})$, and the leader objective values $\Theta_{\I}(\mu,x_{\mu_\Snew}^*)$ and $\Theta_{\I}(\mu,x_{\mu_\W}^*)$ obtained by evaluating the optimistic and pessimistic optimal leader solutions under the different measures.
    \begin{table}[ht]
    \TABLE
    {Comparison of optimal solutions under different follower cooperation behaviors for Example~\ref{ex:behavior-swd}.\label{tab:example-swd-results}}
    {\begin{tabular}{ccccccc}
    \hline
    \textbf{Measure} ($\mu$) & \textbf{Parameters} & \(x^*_{\mu}\) & $\Theta_{\I}(\mu)$ & \(\Theta_{\I}(\mu,x_{\mu_\Snew}^*)\) & \(\Theta_{\I}(\mu,x_{\mu_\W}^*)\) & $\varphi(x^{*}_{\mu})$ \\
    \hline
    {\SBO} -- $\mu_{\Snew}$ & -- & 0.000 & 30.000 & -- & -- & 3.000 \\
    {\WBO} -- $\mu_{\W}$ & -- & 2.000 & 58.000 & -- & -- & 1.000 \\
    {\SWDBO} -- $\mu_{\beta_{\text{Prp}}(\cdot)}$ & -- & 1.600 & 47.200 & 60.000 & 48.000 & 1.400 \\
    {\SWDBO} -- $\mu_{\beta_{\text{Thr}}(\cdot)}$ &  $\delta = 6$ & 1.479 & 46.231 & 58.577 & 48.474 & 1.521 \\
    {\SWDBO} -- $\mu_{\beta_{\text{Str}}(\cdot)}$ & \(a = 5\) & 0.792 & 40.047 & 60.000 & 48.000 & 2.208 \\
    {\SWDBO} -- $\mu_{\beta_{\text{Frg}}(\cdot)}$ & \(a = 0.5\) & 1.750 & 47.675 & 60.000 & 48.000 & 1.250 \\
    {\SWDBO} -- $\mu_{\beta_{\text{Str-p}}(\cdot)}$ & \(p = 1.5\) & 1.343 & 45.207 & 60.000 & 48.000 & 1.657 \\
    {\SWDBO} -- $\mu_{\beta_{\text{Frg-p}}(\cdot)}$ & \(p = 1.5\) & 1.882 & 47.913 & 60.000 & 48.000 & 1.118 \\
    \hline
    \end{tabular}
    }
    {}
    \end{table}
\end{example}

\begin{example}\label{ex:behavior-ibl}
    We extend the problem presented in Example~\ref{ex:behavior-swd} to the {\IBO} general framework, obtaining:
    \begin{equation*}
        \mbox{{\IBO}:} \min_{0 \leq x \leq 2} \left\{ 19x + \int_{\hat{y} \in \mathcal{S}(x,0)} \left( 10\hat{y}_1 + 20\hat{y}_2 \right) \, d\mu_x^\pi (\hat{y}) \right\}, \nonumber
    \end{equation*}
    where ${\mathcal{S}(x,0) = \{y_1 \geq 0,\, y_2 \geq 0 : y_1 + y_2 = 3 - x\}}$, with ${\varphi(x) = 3 - x}$, ${\varphi_{lb} = 1}$, and ${\varphi_{ub} = 3}$, as well as $F_{{\Snew}}(x,\varepsilon) = 10(3-x)$ and $F_{{\W}}(x,\varepsilon)=20(3-x)$ (see Example~\ref{ex:behavior-swd}). All $\hat y \in \mathcal{S}(x,0)$ satisfy $\hat y_2 = 3 - x - \hat y_1$, and thus this problem reduces to:
     \begin{equation*}
        \min_{0 \leq x \leq 2}  \left\{ 60 - x -\int_{\hat{y}_1 \in [0,3-x]} \left(  10\hat{y}_1 \right) \, d\mu_x^\pi (\hat{y}_1)  \right\} = \min_{0 \leq x \leq 2} \left\{  60 - x -\frac{\int_{\hat{y}_1 \in [0,3-x]} \left( 10\hat{y}_1 \right) \pi\!\left( \Tilde{F}(x,\hat{y}_1), \tilde{\varphi}(x)\right) \, d\hat{y}_1 }{ \int_{\hat{y}_1 \in [0,3-x]}  \pi\!\left( \Tilde{F}(x,\hat{y}_1), \tilde{\varphi}(x)\right) \, d\hat{y}_1 } \right\},
    \end{equation*}
    with $\Tilde{F}(x,\hat{y}_1) = -\hat{y}_1/(3-x) + 1/2$ and $\tilde{\varphi}(x) = 1-x/2$. The $\pi(\cdot)$ function in the proprotional case is defined as:
    \begin{align*}
        \pi_{\text{Prp}}\!\left( \Tilde{F}(x,\hat{y}), \tilde{\varphi}(x) \right) &= \pi_c + \left(1-x \right) \left(-\tfrac{\hat{y}_1}{3-x} + \tfrac{1}{2} \right), 
    \end{align*}
    As in Example~\ref{ex:behavior-swd}, Table~\ref{tab:example-ibl-results} reports the leader’s optimal decision $x_{\mu}^*$ and optimal objective $\Theta_{\I}(\mu)$ for each considered measure. It shows the objective values $\Theta_{\I}(\mu,x_{\mu_\Snew}^*)$ and $\Theta_{\I}(\mu,x_{\mu_\W}^*)$ obtained by evaluating the optimistic and pessimistic leader solutions under the different follower-cooperation approaches, as well as the optimal value of the follower $\varphi(x_{\mu}^*)$.
    \begin{table}[ht]
    \TABLE
    {Comparison of optimal solutions under different follower cooperation behaviors for Example~\ref{ex:behavior-ibl}.\label{tab:example-ibl-results}}
    {\begin{tabular}{ccccccc}
    \hline
    \textbf{Measure} ($\mu$) & \textbf{Parameters} & \(x^*_{\mu}\) & $\Theta_{\I}(\mu)$ & \(\Theta_{\I}(\mu,x_{\mu_\Snew}^*)\) & \(\Theta_{\I}(\mu,x_{\mu_\W}^*)\) & $\varphi(x^{*}_{\mu})$ \\
    \hline
    {\IBO} -- $\mu_{\pi_{\text{Ntr}}(\cdot)}$ & -- & 0.000 & 45.000 & 45.000 & 53.000 & 3.000 \\
    {\IBO} -- $\mu_{\pi_{\text{Prp}}(\cdot)}$ & & 0.800 & 48.934 & 50.000 & 51.334 & 2.200 \\
    \hline
    \end{tabular}
    }
    {}
    \end{table}
\end{example}

\section{Detailed {\MCMC} Hit-and-Run Procedure}\label{app:details-har}

In this section, we describe the {\MCMC} Hit-and-Run procedure, presented in Algorithm~\ref{alg:hitrun} and used to sample from $\mathcal{S}_{\text{Lin}}(x,\varepsilon)$ following the generalized measures $\mu = \mu_{\pi(\cdot)}$. In this case, the exogenous random vector is \({\boldsymbol{\zeta} = \big( \boldsymbol{\zeta}_{\omega_1} = (\boldsymbol{\upsilon}_{\omega_1},\boldsymbol{\tau}_{\omega_1}), \ldots, \boldsymbol{\zeta}_{\omega_N} =(\boldsymbol{\upsilon}_{\omega_N},\boldsymbol{\tau}_{\omega_N})\big)}\) with realization $\zeta$, whose components are independent and identically distributed. The direction random component satisfies either $\boldsymbol{\upsilon}_\omega \sim U(\mathbb{S}^{n_y-1})$ (standard Hit-and-Run), where $\mathbb{S}^{n_y-1}$ denotes the unit Euclidean sphere in $\mathbb{R}^{n_y}$, or $\boldsymbol{\upsilon}_\omega \sim U(\{u_1,\dots,u_{n_y}\})$ (coordinate Hit-and-Run), where $u_\ell$ denotes the $\ell$-th canonical basis vector in $\mathbb{R}^{n_y}$. Moreover, the step size random component is $\boldsymbol{\tau}_\omega \sim U(0,1)$. Finally, if $\mathcal{S}_{\mathrm{Lin}}(x,\varepsilon)$ includes equality constraints, the direction $\boldsymbol{\upsilon}_\omega$ is sampled in the null space of the equality system to preserve feasibility. In our experimental results, we use the coordinate Hit-and-Run.

\begin{algorithm}[htb]
\caption{{\MCMC} Hit-and-Run} 
\label{alg:hitrun}
\small
\Input{$N$, $\{\zeta_{\omega}\}_{\omega \in \mathcal{W}_N} = \left\{\left(\upsilon_{\omega}, \tau_{\omega} \right)\right\}_{\omega \in \mathcal{W}_N}$, $\varepsilon$, $\gamma$, $K$, $L$, $\pi(\cdot) \in \Pi$, $\mu = \mu_{\pi(\cdot)}$, $x \in \mathcal{X}$, $\hat{y}_{\texttt{C}} \in \mathcal{S}_{\mathrm{Lin}}(x,\varepsilon)$, $\tilde\varphi(x)$;}

\Output{$\{\hat{y}_{\omega}\}_{\omega \in \mathcal{W}_N}$;}

\textbf{Step 0:} Initialize $\hat{y}_0 \gets \hat{y}_{\texttt{C}}$\;

\For{$i = 1, \ldots, N$}{

    \textbf{Step 1:} Compute the feasible interval $[\alpha_{\min}, \alpha_{\max}] = \{\alpha \in [-L,L] : \hat{y}_{\omega_{i-1}} + \alpha \upsilon_{\omega_{i}} \in \mathcal{S}_{\mathrm{Lin}}(x,\varepsilon)\}$\;

    \textbf{Step 2:} Compute $\alpha_{\omega_i} = \mathcal{P}_{\pi(\cdot)}^{-1}(\tau_{\omega_i})$, where  $\mathcal{P}_{\pi(\cdot)}(\alpha)$ is the cdf of $\pi \left( \tilde F(x,\hat{y}_{\omega_{i-1}} + \alpha \upsilon_{\omega_i}), \tilde\varphi(x) \right)$ with parameter $\gamma$ along the segment $[\alpha_{\min}, \alpha_{\max}]$\;

    \textbf{Step 3:} Set $\hat{y}_{\omega_i} \gets \hat{y}_{\omega_{i-1}} + \alpha_{\omega_i} \upsilon_{\omega_i}$\;
}
\end{algorithm}

For a leader decision $x \in \mathcal{X}$, the Markov chain $\{\hat y_{\omega_i}\}_{\omega \in \mathcal{W}_N}$ on $\mathcal{S}_{\mathrm{Lin}}(x,\varepsilon)$ is defined recursively by
\begin{equation}\label{eq:markov-chain}
\hat y_0 = \hat y_{\texttt{C}}, \qquad \hat y_{\omega_i} = mc(\hat y_{\omega_{i-1}}, x, \zeta_{\omega_i}), \quad i=1,\ldots,N, 
\end{equation}
where $\hat y_{\texttt{C}}$ is the Chebyshev center of $\mathcal{S}_{\mathrm{Lin}}(x,\varepsilon)$ and $mc (\cdot)$ is the Hit-and-Run transition mapping. This mapping induces an endogenous transition kernel $K_x(\hat y' \mid \hat y)$ with stationary measure $\mu_x$, that is,
\begin{equation}
    \mu_x(A) = \int_{\hat y \in \mathcal{S}(x,\varepsilon)} K_x(\hat y' \in A \mid \hat y) \, d\mu_x(\hat y), \qquad \forall A \subseteq \mathcal{S}_{\mathrm{Lin}}(x,\varepsilon). \nonumber
\end{equation}
The transformation \(t_{\mu}(x,\zeta) = (\hat y_{\omega_1},\ldots,\hat y_{\omega_N}) \) is therefore implicitly defined by the Markov recursion~\eqref{eq:markov-chain} for the realization $\zeta$ of $\boldsymbol{\zeta}$. Although the sequence of samples $\{\hat y_{\omega}\}_{\omega \in \mathcal{W}_N}$ is dependent, the chain's irreducibility and aperiodicity ensure ergodicity; thus, empirical averages computed from $t_{\mu}(x,\zeta)$ converge to expectations under $\mu_x$ as $N$ increases. The proposal kernel is symmetric by construction, due to the sampling of the direction $\upsilon_{\omega_i}$ and the definition of the step $\alpha_{\omega_i}$. The target measure $\mu_x$ is thus proportional to $\pi(\cdot)$.

In Step~1 at iteration $i$ of Algorithm~\ref{alg:hitrun}, the bounds $[\alpha_{\text{min}},\alpha_{\text{max}}]$ are defined as:
\begin{align*}
    &\alpha_{\text{min}} = \max \{ \alpha : \, \hat{y}_{\omega_{i-1}} + \alpha \upsilon_{\omega_i} \in \mathcal{S}_{\text{Lin}}(x,\varepsilon), \, -L \leq \alpha \leq 0 \} , \\
    &\alpha_{\text{max}} = \min \{ \alpha : \, \hat{y}_{\omega_{i-1}} + \alpha \upsilon_{\omega_i} \in \mathcal{S}_{\text{Lin}}(x,\varepsilon), \, 0 \leq \alpha \leq L  \}.
\end{align*}
These bounds represent the largest feasible steps in the direction $\upsilon_{\omega_i}$ from the current point $\hat{y}$ that remain within the set $\mathcal{S}_{\text{Lin}}(x,\varepsilon)$, while the constants $-L$ and $L$ provide explicit global bounds on the step length. We note that these step-size bounds can be written explicitly as
\begin{align}
    &\alpha_{\text{min}}
   = \max \Bigg\{ -L, 
       \max_{ \substack{j=1,\ldots,n_y\\ \upsilon_{\omega_i}^{j}>0}}
           \left( -\frac{\hat y_{\omega_{i-1}}^j}{\upsilon_{\omega_i}^{j}} \right),
       \;
       \max_{\substack{l=1,\ldots,m_2\\ D_l \upsilon_{\omega_i} < 0}}
           \left( \frac{(a - Cx)_l - D_l \hat y_{\omega_i}}{D_l \upsilon_{\omega_{i-1}}} \right)
     \Bigg\},\\
    &\alpha_{\text{max}}
   = \min \Bigg\{ L, 
       \min_{\substack{j=1,\ldots,n_y\\ \upsilon_{\omega_i}^{j}<0}}
           \left( -\frac{\hat y_{\omega_{i-1}}^j}{\upsilon_{\omega_{i}}^{j}} \right),
       \;
       \min_{\substack{l=1,\ldots,L\\ D_l \upsilon_{\omega_i} > 0}}
           \left( \frac{(a - Cx)_l - D_l \hat y_{\omega_{i-1}}}{D_l \upsilon_{\omega_i}} \right)
     \Bigg\}.
\end{align}
Furthermore, in Step~2 of the algorithm, the step $\alpha_{\omega_i}$ along the line segment is sampled proportionally to $\pi(\cdot) \in \Pi$ using the inverse cdf method, as a closed-form inverse is available in all considered cases. When such a closed form is unavailable or difficult to compute, alternative sampling techniques can be used. 

A rounding procedure can be applied to rescale the set $\mathcal{S}_{\mathrm{Lin}}(x,\varepsilon)$ so that its width is approximately uniform across all dimensions, improving the mixing or convergence of the Hit-and-Run sampler~\citep{lovasz2006hit}. Algorithm~\ref{alg:hitrun} illustrates the sampling of a single Markov chain; however, multiple chains initialized at the Chebyshev center can be employed. Thinning can also be applied to reduce autocorrelation among the samples generated by the Markov chain by retaining only every few draws. For fixed $x \in \mathcal{X}$, if $F(x,\boldsymbol{\hat{y}})$ is linear in $\boldsymbol{\hat{y}}$, then $\pi(\cdot) \in \Pi$ is log-concave in $\boldsymbol{\hat{y}}$. In this case, the mixing time of Hit-and-Run is polynomial in the dimension $n_y$~\citep{lovasz2006hit}, ensuring efficient convergence. Moreover, we pre-sample $\{\zeta_\omega\}_{\omega \in \mathcal{W}_N}$ and reuse them whenever a feasible leader solution is evaluated. This induces correlation between estimators $\hat{g}^{\mu}_{N}(x)$ across different leader decisions $x$, reducing the variance of their differences and improving approximation stability and efficiency~\citep{bazotte2025solving}.

Once the chain has mixed or converged, its subsequent states can be regarded as approximately drawn from the stationary distribution, enabling the use of standard {\MCMC} techniques within the {\SAA} program. Accordingly, we define a single Markov chain for the {\SAA} program, whose states correspond to scenarios $\omega \in \mathcal{W}_N$, each characterized by $\zeta_\omega = (\upsilon_\omega,\tau_\omega,\phi_\omega)$.

\section{Detailed Instances}~\label{app:detailed-instances}

We adapt 170 instances from two groups of the BOBILib library~\citep{BOBILib:2026}. The first group consists of 110 instances from the \texttt{denegre} set~\citep{denegre2011interdiction}, with no leader constraints and either 20 or 30 follower constraints. The leader has $n_x \in \{5,15,20,25\}$ integer variables, taking values in $\{0,\dots,1500\}$ (equivalent to $11\,n_x$ binary variables) or in $\{0,\dots,300\}$ (equivalent to $9\,n_x$ binary variables). The follower has $n_y \in \{5,10,15\}$ variables, relaxed to be continuous in $[0,1500]$ or $[0,300]$. The second group comprises 60 instances from the \texttt{xuwang} set~\citep{xu2014exact}, with ${n_x \!= \!n_y \!\in\! \{10,60,110,160,210,260\}}$, integer leader variables in $\{0,\dots,10\}$ (equivalent to $4\, n_x$ binary variables) and follower variables relaxed to be continuous and nonnegative. The number of constraints at each level is $4, 24, 44, 64, 84,$ and $104$, where the initial upper-level coupling constraints are modified \emph{a priori} to remove their dependence on follower variables. Note that we use general-purpose benchmark instances rather than application-specific datasets to evaluate the proposed framework in a broad and unbiased setting. This allows us to assess its performance independently of particular application structures and verify its general applicability. For all instances, the follower’s objective coefficients are set to zero with probability 1/2, otherwise the original value is kept. The aim for this procedure is to increase the likelihood of multiple optimal follower solutions.

\section{Detailed Experimental Results}~\label{app:additional-results}

In this section, we present the experimental results. The analysis is organized into two parts: (i) computational performance (Section~\ref{app:additional-results-perf}) and (ii) comparison of follower tie-breaking behaviors (Section~\ref{app:additional-results-comp}).

\subsection{Computational performance}~\label{app:additional-results-perf}

In this section, we present detailed experimental results on the computational performance of the proposed methods for the different measures. Specifically, we report the optimality gap [$\text{gap}_{opt}$($\%$)] (as defined in Section~\ref{subsubsec:baseline-performance}) and runtime of the {\DE} program for strong-weak fixed measures (Table~\ref{tab:details-res-sw}); the {\SAA} program (with $N=100$ scenarios) and the {\DE} program for strong-weak decision-dependent measures (Tables~\ref{tab:details-res-saa-swd} and~\ref{tab:details-res-dep-swd}); and the {\SAA} program (with $N=1500$ scenarios generated from a single Markov chain and a thinning factor of ten) for generalized decision-dependent measures (Table~\ref{tab:details-res-gen}), each evaluated per number of leader and follower variables. We also report the runtime of the {\MCMC} Hit-and-Run method for $N'=4\cdot10^4$ samples, generated using four Markov chains with a thinning factor of 50, also per number of leader and follower variables (Table~\ref{tab:details-res-gen-eval}). Finally, Figure~\ref{fig:comp-de-sw} shows the runtime (Figure~\ref{subfig:comp-de-sw-time}) and optimality gap (Figure~\ref{subfig:swd-saa-de-gap}) of the {\DE} program for strong-weak decision-dependent measures per percentage of instances.
Note that, for the proposed instances, each $x_j$ can be represented using ${s_j = \lceil \log_2(u_j - l_j + 1) \rceil}$ binary variables $x'_0, \dots, x'_{h_j-1}$, with $x_j = l_j + x'_0 + 2 x'_1 + 4 x'_2 + \ldots + 2^{h_j-1} x'_{{s_j}-1}$, and $l_j \le x_j \le u_j$ for ${j = 1, \dots, n_x}$.

\begin{table}[]
\TABLE
{Detailed results of the {\DE} program for strong-weak fixed measures. \label{tab:details-res-sw}}
{
\begin{adjustbox}{width=\textwidth}
\begin{tabular}{cc|cc|cc|cc|cc|cc}
\hline
\multirow{2}{*}{$n_x$} & \multicolumn{1}{l|}{\multirow{2}{*}{$n_y$}} & \multicolumn{2}{c|}{$\beta^1 = \mu_{\Snew}$} & \multicolumn{2}{c|}{$\beta^{0.7}$}  & \multicolumn{2}{c|}{$\beta^{0.5}$}  & \multicolumn{2}{c|}{$\beta^{0.3}$} & \multicolumn{2}{c}{$\beta^{0}=\mu_{\W}$} \\ \cline{3-12} 
                       & \multicolumn{1}{l|}{}                       & $\text{gap}_{opt}$($\%$)      & Time (s)            & $\text{gap}_{opt}$($\%$) & Time (s)        & $\text{gap}_{opt}$($\%$) & Time (s)        & $\text{gap}_{opt}$($\%$)  & Time (s)      & $\text{gap}_{opt}$($\%$)     & Time (s)         \\ \hline
5                      & 10                     & 0.0$\pm$0.0                             & 0.1$\pm$0.2                     & 0.0$\pm$0.0                             & 0.1$\pm$0.2                     & 0.0$\pm$0.0                             & 0.0$\pm$0.0                     & 0.0$\pm$0.0                             & 0.0$\pm$0.0                     & 0.0$\pm$0.0                             & 0.0$\pm$0.0                     \\
5                      & 15                     & 0.0$\pm$0.0                             & 0.4$\pm$0.5                     & 0.0$\pm$0.0                             & 1.1$\pm$1.6                     & 0.0$\pm$0.0                             & 0.7$\pm$0.9                     & 0.0$\pm$0.0                             & 0.4$\pm$0.4                     & 0.0$\pm$0.0                             & 0.2$\pm$0.3                     \\
10                     & 10                     & 0.0$\pm$0.0                             & 0.5$\pm$0.5                     & 0.0$\pm$0.0                             & 1.8$\pm$2.2                     & 0.0$\pm$0.0                             & 0.1$\pm$0.2                     & 0.0$\pm$0.0                             & 0.1$\pm$0.2                     & 0.0$\pm$0.0                             & 0.0$\pm$0.0                     \\
15                     & 5                      & 0.0$\pm$0.0                             & 0.0$\pm$0.0                     & 0.0$\pm$0.0                             & 0.1$\pm$0.1                     & 0.0$\pm$0.0                             & 0.0$\pm$0.1                     & 0.0$\pm$0.0                             & 0.1$\pm$0.1                     & 0.0$\pm$0.0                             & 0.0$\pm$0.1                     \\
20                     & 10                     & 0.0$\pm$0.0                             & 0.6$\pm$0.3                     & 0.0$\pm$0.0                             & 1.3$\pm$0.8                     & 0.0$\pm$0.0                             & 1.4$\pm$1.0                     & 0.0$\pm$0.0                             & 1.5$\pm$0.9                     & 0.0$\pm$0.0                             & 1.2$\pm$0.9                     \\
25                     & 5                      & 0.0$\pm$0.0                             & 46.0$\pm$3.8                    & 0.0$\pm$0.0                             & 91.8$\pm$6.7                    & 0.0$\pm$0.0                             & 98.4$\pm$6.6                    & 0.0$\pm$0.0                             & 105.2$\pm$7.5                   & 0.0$\pm$0.0                             & 104.4$\pm$7.9                   \\
10                     & 10                     & 0.0$\pm$0.0                             & 0.0$\pm$0.0                     & 0.0$\pm$0.0                             & 0.0$\pm$0.0                     & 0.0$\pm$0.0                             & 0.0$\pm$0.0                     & 0.0$\pm$0.0                             & 0.0$\pm$0.0                     & 0.0$\pm$0.0                             & 0.0$\pm$0.0                     \\
60                     & 60                     & 0.0$\pm$0.0                             & 0.9$\pm$0.6                     & 0.0$\pm$0.0                             & 1.6$\pm$0.8                     & 0.0$\pm$0.0                             & 1.4$\pm$1.0                     & 0.0$\pm$0.0                             & 0.6$\pm$0.4                     & 0.0$\pm$0.0                             & 0.8$\pm$0.5                     \\
110                    & 110                    & 0.0$\pm$0.0                             & 10.8$\pm$4.7                    & 0.0$\pm$0.0                             & 13.1$\pm$9.1                    & 0.0$\pm$0.0                             & 6.2$\pm$3.6                     & 0.0$\pm$0.0                             & 3.4$\pm$2.4                     & 0.0$\pm$0.0                             & 2.8$\pm$1.9                     \\
160                    & 160                    & 0.0$\pm$0.0                             & 39.4$\pm$25.4                   & 0.0$\pm$0.0                             & 42.4$\pm$39.9                   & 0.0$\pm$0.0                             & 19.0$\pm$13.1                   & 0.0$\pm$0.0                             & 13.1$\pm$8.3                    & 0.0$\pm$0.0                             & 7.1$\pm$4.5                     \\
210                    & 210                    & 0.0$\pm$0.0                             & 447.3$\pm$289.7                 & 0.0$\pm$0.0                             & 146.2$\pm$105.1                 & 0.0$\pm$0.0                             & 30.2$\pm$22.6                   & 0.0$\pm$0.0                             & 22.7$\pm$22.8                   & 0.0$\pm$0.0                             & 7.6$\pm$3.3                     \\
260                    & 260                    & 0.0$\pm$0.0                             & 2191.8$\pm$2194.0               & 0.0$\pm$0.0                             & 1140.0$\pm$786.9                & 0.0$\pm$0.0                             & 349.9$\pm$283.3                 & 0.0$\pm$0.0                             & 117.8$\pm$123.2                 & 0.0$\pm$0.0                             & 55.0$\pm$68.1                   \\ \hline
\end{tabular}
\end{adjustbox}
}
{The values are averages over instance-behavior configurations with the same variable size, with a 95\% confidence interval.}
\end{table}

\begin{table}[]
\TABLE
{Detailed results of the {\SAA} program for strong-weak decision-dependent measures.\label{tab:details-res-saa-swd}}
{
\begin{adjustbox}{width=\textwidth}
\centering
\begin{tabular}{cc|cc|cc|cc|cc|cc|cc}
\hline
\multirow{2}{*}{$n_x$} & \multirow{2}{*}{$n_y$} & \multicolumn{2}{c|}{$\beta_{\text{Prp}}(\cdot)$} & \multicolumn{2}{c|}{$\beta_{\text{Thr}}(\cdot)$} & \multicolumn{2}{c|}{$\beta_{\text{Str}}(\cdot)$} & \multicolumn{2}{c|}{$\beta_{\text{Frg}}(\cdot)$} & \multicolumn{2}{c|}{$\beta_{\text{Str-p}}(\cdot)$} & \multicolumn{2}{c}{$\beta_{\text{Frg-p}}(\cdot)$} \\ \cline{3-14} 
                       &                        & $\text{gap}_{opt}$($\%$)        & Time (s)              & $\text{gap}_{opt}$($\%$)        & Time (s)              & $\text{gap}_{opt}$($\%$)    & Time (s)           & $\text{gap}_{opt}$($\%$)    & Time (s)           & $\text{gap}_{opt}$($\%$)    & Time (s)           & $\text{gap}_{opt}$($\%$)    & Time (s)          \\ \hline
5                      & 10                     & 0.0$\pm$0.0                             & 1.1$\pm$0.6                     & 0.0$\pm$0.0                             & 0.6$\pm$0.2                     & 0.0$\pm$0.0                             & 3.2$\pm$1.7                     & 0.0$\pm$0.0                             & 0.9$\pm$0.2                     & 0.0$\pm$0.0                             & 4.1$\pm$2.1                     & 0.0$\pm$0.0 & 0.6$\pm$0.2     \\
5                      & 15                     & 0.0$\pm$0.0                             & 4.0$\pm$3.4                     & 0.0$\pm$0.0                             & 4.0$\pm$2.3                     & 0.0$\pm$0.0                             & 57.8$\pm$56.2                   & 0.0$\pm$0.0                             & 2.0$\pm$0.7                     & 0.0$\pm$0.0                             & 114.7$\pm$120.9                 & 0.0$\pm$0.0 & 1.5$\pm$0.4     \\
10                     & 10                     & 0.0$\pm$0.0                             & 2.1$\pm$2.5                     & 0.0$\pm$0.0                             & 2.2$\pm$2.1                     & 0.0$\pm$0.0                             & 30.2$\pm$30.9                   & 0.0$\pm$0.0                             & 1.1$\pm$0.4                     & 0.0$\pm$0.0                             & 45.2$\pm$44.3                   & 0.0$\pm$0.0 & 0.8$\pm$0.3     \\
15                     & 5                      & 0.0$\pm$0.0                             & 1.2$\pm$0.8                     & 0.0$\pm$0.0                             & 1.3$\pm$1.2                     & 0.0$\pm$0.0                             & 1.3$\pm$0.4                     & 0.0$\pm$0.0                             & 0.8$\pm$0.2                     & 0.0$\pm$0.0                             & 1.4$\pm$0.6                     & 0.0$\pm$0.0 & 0.5$\pm$0.2     \\
20                     & 10                     & 0.0$\pm$0.0                             & 17.9$\pm$7.1                    & 0.0$\pm$0.0                             & 15.0$\pm$4.2                    & 0.0$\pm$0.0                             & 19.0$\pm$4.4                    & 0.0$\pm$0.0                             & 13.4$\pm$2.2                    & 0.0$\pm$0.0                             & 28.2$\pm$15.4                   & 0.0$\pm$0.0 & 13.5$\pm$3.1    \\
25                     & 5                      & 0.0$\pm$0.0                             & 241.5$\pm$131.8                 & 0.0$\pm$0.0                             & 192.3$\pm$77.6                  & 0.0$\pm$0.0                             & 2015.5$\pm$464.3                & 0.0$\pm$0.0                             & 137.6$\pm$18.8                  & 0.0$\pm$0.0                             & 2286.3$\pm$525.7                & 0.0$\pm$0.0 & 127.7$\pm$10.2  \\
10                     & 10                     & 0.0$\pm$0.0                             & 0.6$\pm$0.4                     & 0.0$\pm$0.0                             & 0.4$\pm$0.2                     & 0.0$\pm$0.0                             & 0.7$\pm$0.2                     & 0.0$\pm$0.0                             & 0.8$\pm$0.2                     & 0.0$\pm$0.0                             & 0.6$\pm$0.2                     & 0.0$\pm$0.0 & 0.4$\pm$0.2     \\
60                     & 60                     & 0.0$\pm$0.0                             & 7.3$\pm$3.7                     & 0.0$\pm$0.0                             & 5.5$\pm$1.9                     & 0.0$\pm$0.0                             & 28.4$\pm$12.7                   & 0.0$\pm$0.0                             & 4.8$\pm$1.0                     & 0.0$\pm$0.0                             & 40.0$\pm$18.2                   & 0.0$\pm$0.0 & 4.1$\pm$1.1     \\
110                    & 110                    & 0.0$\pm$0.0                             & 40.5$\pm$41.6                   & 0.0$\pm$0.0                             & 33.9$\pm$20.1                   & 0.0$\pm$0.0                             & 468.5$\pm$477.8                 & 0.0$\pm$0.0                             & 25.6$\pm$13.9                   & 0.2$\pm$0.4                             & 884.4$\pm$814.4                 & 0.0$\pm$0.0 & 13.5$\pm$6.2    \\
160                    & 160                    & 0.0$\pm$0.0                             & 179.1$\pm$156.6                 & 0.0$\pm$0.0                             & 254.6$\pm$209.5                 & 2.2$\pm$2.8                             & 1860.3$\pm$1071.1               & 0.0$\pm$0.0                             & 78.1$\pm$35.3                   & 4.2$\pm$5.9                             & 2883.6$\pm$1400.2               & 0.0$\pm$0.0 & 53.4$\pm$30.2   \\
210                    & 210                    & 0.0$\pm$0.0                             & 412.9$\pm$362.2                 & 0.1$\pm$0.1                             & 630.1$\pm$590.0                 & 4.9$\pm$2.6                             & 4899.1$\pm$1549.8               & 0.0$\pm$0.0                             & 140.1$\pm$82.1                  & 10.4$\pm$7.1                            & 7165.8$\pm$1714.9               & 0.0$\pm$0.0 & 72.4$\pm$60.1   \\
260                    & 260                    & 1.6$\pm$3.7                             & 3379.4$\pm$2991.3               & 1.1$\pm$1.1                             & 2689.3$\pm$1183.8               & 13.7$\pm$6.6                            & 6639.9$\pm$1524.7               & 0.0$\pm$0.0                             & 1292.6$\pm$780.4                & 26.7$\pm$12.6                           & 8466.2$\pm$1617.4               & 0.0$\pm$0.0 & 494.9$\pm$403.9 \\ \hline
\end{tabular}
\end{adjustbox}
}
{The values are averages over instance-behavior configurations with the same variable size, with a 95\% confidence interval.}
\end{table}

\begin{table}[]
\TABLE
{Detailed results of the {\DE} program for strong-weak decision-dependent measures.\label{tab:details-res-dep-swd}}
{
\begin{adjustbox}{width=\textwidth}
\begin{tabular}{cc|cc|cc|cc|cc|cc|cc}
\hline
\multirow{2}{*}{$n_x$} & \multirow{2}{*}{$n_y$} & \multicolumn{2}{c|}{$\beta_{\text{Prp}}(\cdot)$} & \multicolumn{2}{c|}{$\beta_{\text{Thr}}(\cdot)$} & \multicolumn{2}{c|}{$\beta_{\text{Str}} (\cdot)$} & \multicolumn{2}{c|}{$\beta_{\text{Frg}} (\cdot)$} & \multicolumn{2}{c|}{$\beta_{\text{Str-p}} (\cdot)$} & \multicolumn{2}{c}{$\beta_{\text{Frg-p}} (\cdot)$} \\ \cline{3-14} 
                       &                        & $gap_{opt}$($\%$)        & Time (s)              & $gap_{opt}$($\%$)        & Time (s)              & $gap_{opt}$($\%$)    & Time (s)           & $gap_{opt}$($\%$)    & Time (s)           & $gap_{opt}$($\%$)     & Time (s)            & $gap_{opt}$($\%$)     & Time (s)           \\ \hline
5                      & 10                     & 0.0$\pm$0.0                             & 0.6$\pm$0.5                     & 0.0$\pm$0.0                             & 18.3$\pm$2.3                    & 0.0$\pm$0.0                             & 34.4$\pm$18.3                   & 0.0$\pm$0.0                             & 16.6$\pm$1.0                    & 0.0$\pm$0.0                             & 24.8$\pm$6.6                    & 0.0$\pm$0.0 & 15.7$\pm$0.6      \\
5                      & 15                     & 0.0$\pm$0.0                             & 9.4$\pm$11.7                    & 0.0$\pm$0.0                             & 34.2$\pm$16.9                   & 0.0$\pm$0.0                             & 98.8$\pm$63.2                   & 0.0$\pm$0.0                             & 20.5$\pm$2.4                    & 0.0$\pm$0.0                             & 209.0$\pm$151.6                 & 0.0$\pm$0.0 & 19.0$\pm$1.9      \\
10                     & 10                     & 0.0$\pm$0.0                             & 4.8$\pm$6.5                     & 0.0$\pm$0.0                             & 27.2$\pm$11.2                   & 0.0$\pm$0.0                             & 132.1$\pm$153.3                 & 0.0$\pm$0.0                             & 19.5$\pm$2.7                    & 0.0$\pm$0.0                             & 148.0$\pm$132.6                 & 0.0$\pm$0.0 & 24.1$\pm$10.5     \\
15                     & 5                      & 0.0$\pm$0.0                             & 1.0$\pm$1.3                     & 0.0$\pm$0.0                             & 23.1$\pm$12.9                   & 0.0$\pm$0.0                             & 21.4$\pm$5.5                    & 0.0$\pm$0.0                             & 16.1$\pm$0.7                    & 0.0$\pm$0.0                             & 24.3$\pm$11.0                   & 0.0$\pm$0.0 & 15.8$\pm$0.8      \\
20                     & 10                     & 0.0$\pm$0.0                             & 29.6$\pm$21.0                   & 0.0$\pm$0.0                             & 67.1$\pm$36.0                   & 0.0$\pm$0.0                             & 46.6$\pm$5.9                    & 0.0$\pm$0.0                             & 46.7$\pm$5.8                    & 0.0$\pm$0.0                             & 64.4$\pm$38.7                   & 0.0$\pm$0.0 & 40.9$\pm$5.1      \\
25                     & 5                      & 0.0$\pm$0.0                             & 2321.9$\pm$1024.9               & 0.0$\pm$0.0                             & 1654.4$\pm$348.5                & 0.0$\pm$0.0                             & 3514.2$\pm$588.7                & 0.0$\pm$0.0                             & 1258.6$\pm$314.8                & 0.0$\pm$0.0                             & 3570.7$\pm$639.0                & 0.0$\pm$0.0 & 1777.3$\pm$512.1  \\
10                     & 10                     & 0.0$\pm$0.0                             & 0.0$\pm$0.0                     & 0.0$\pm$0.0                             & 14.4$\pm$0.3                    & 0.0$\pm$0.0                             & 14.5$\pm$0.3                    & 0.0$\pm$0.0                             & 14.2$\pm$0.4                    & 0.0$\pm$0.0                             & 14.0$\pm$0.4                    & 0.0$\pm$0.0 & 13.9$\pm$0.4      \\
60                     & 60                     & 0.0$\pm$0.0                             & 13.4$\pm$6.9                    & 0.0$\pm$0.0                             & 35.7$\pm$8.8                    & 0.0$\pm$0.0                             & 111.5$\pm$48.6                  & 0.0$\pm$0.0                             & 34.0$\pm$6.4                    & 0.0$\pm$0.0                             & 145.1$\pm$80.2                  & 0.0$\pm$0.0 & 29.7$\pm$10.6     \\
110                    & 110                    & 0.0$\pm$0.0                             & 97.9$\pm$73.1                   & 0.0$\pm$0.0                             & 137.6$\pm$49.8                  & 0.6$\pm$0.9                             & 1494.8$\pm$956.6                & 0.0$\pm$0.0                             & 104.8$\pm$43.7                  & 2.0$\pm$2.5                             & 2024.0$\pm$1222.6               & 0.0$\pm$0.0 & 76.4$\pm$34.7     \\
160                    & 160                    & 0.0$\pm$0.0                             & 226.9$\pm$172.1                 & 0.0$\pm$0.0                             & 688.8$\pm$346.8                 & 6.7$\pm$5.9                             & 3839.6$\pm$1429.9               & 0.0$\pm$0.0                             & 519.9$\pm$202.1                 & 8.3$\pm$6.0                             & 5317.3$\pm$1657.8               & 0.0$\pm$0.0 & 405.7$\pm$229.7   \\
210                    & 210                    & 0.0$\pm$0.0                             & 703.1$\pm$494.1                 & 0.3$\pm$0.4                             & 1159.2$\pm$754.5                & 10.4$\pm$5.4                            & 6221.4$\pm$1545.5               & 0.0$\pm$0.0                             & 644.4$\pm$254.7                 & 21.6$\pm$8.3                            & 8300.4$\pm$1584.9               & 0.0$\pm$0.0 & 467.3$\pm$365.7   \\
260                    & 260                    & 5.2$\pm$8.0                             & 4817.1$\pm$3276.3               & 7.6$\pm$7.2                             & 3948.0$\pm$1328.1               & 23.6$\pm$9.9                            & 7786.8$\pm$1361.1               & 2.0$\pm$1.6                             & 3434.2$\pm$1308.3               & 33.2$\pm$10.9                           & 8855.5$\pm$1487.1               & 1.3$\pm$1.5 & 2302.4$\pm$1292.6 \\ \hline
\end{tabular}
\end{adjustbox}
}
{The values are averages over instance-behavior configurations with the same variable size, with a 95\% confidence interval.}
\end{table}

\begin{table}
\TABLE
{Detailed results for generalized decision-dependent measures.\label{tab:details-res-gen}}
{
\begin{adjustbox}{width=\textwidth}
\begin{tabular}{cc|cc|cc|cc|cc|cc}
\hline
\multirow{2}{*}{$n_x$} & \multirow{2}{*}{$n_y$} & \multicolumn{2}{c|}{$\pi_{\text{Ntr}}(\cdot)$} & \multicolumn{2}{c|}{$\pi_{\text{Prp}}(\cdot)$} & \multicolumn{2}{c|}{$\pi_{\text{Frg}}(\cdot)$} & \multicolumn{2}{c|}{$\pi_{\text{Frg}^+}(\cdot)$} & \multicolumn{2}{c}{$\pi_{\text{Str}}(\cdot)$} \\ \cline{3-12} 
                       &                        & $\text{gap}_{opt}$($\%$)       & Time (s)             & $\text{gap}_{opt}$($\%$)       & Time (s)             & $\text{gap}_{opt}$($\%$)    & Time (s)         & $\text{gap}_{opt}$($\%$)    & Time (s)           & $\text{gap}_{opt}$($\%$)   & Time (s)          \\ \hline
55                     & 10                     & 2.1$\pm$4.8                             & 1144.7$\pm$2428.0               & 1.7$\pm$3.7                             & 1143.0$\pm$2428.3               & 0.9$\pm$1.0                             & 530.2$\pm$522.5                 & 0.4$\pm$0.9                             & 296.4$\pm$432.2                 & 1.9$\pm$1.7                             & 1138.5$\pm$924.9                \\
55                     & 15                     & 0.0$\pm$0.0                             & 1711.6$\pm$1737.9               & 0.0$\pm$0.0                             & 1739.3$\pm$1758.1               & 0.0$\pm$0.0                             & 1867.3$\pm$590.0                & 0.0$\pm$0.0                             & 1894.4$\pm$683.2                & 0.0$\pm$0.0                             & 1776.3$\pm$722.2                \\
110                    & 10                     & 13.8$\pm$14.0                           & 4691.0$\pm$3775.6               & 20.6$\pm$23.3                           & 4767.5$\pm$3734.0               & 15.0$\pm$5.4                            & 4684.6$\pm$1203.2               & 14.5$\pm$6.0                            & 4704.8$\pm$1433.7               & 14.4$\pm$5.8                            & 4913.8$\pm$1445.2               \\
165                    & 5                      & 4.2$\pm$4.0                             & 3133.6$\pm$2225.5               & 4.7$\pm$4.6                             & 2884.6$\pm$2206.4               & 4.2$\pm$1.4                             & 2886.9$\pm$773.6                & 4.4$\pm$1.7                             & 2801.2$\pm$923.7                & 4.2$\pm$1.7                             & 2988.5$\pm$919.8                \\
180                    & 10                     & 0.0$\pm$0.0                             & 9.6$\pm$1.2                     & 0.0$\pm$0.0                             & 9.5$\pm$1.2                     & 0.0$\pm$0.0                             & 9.3$\pm$0.4                     & 0.0$\pm$0.0                             & 9.6$\pm$0.6                     & 0.0$\pm$0.0                             & 9.4$\pm$0.5                     \\
225                    & 5                      & 0.0$\pm$0.0                             & 130.5$\pm$14.5                  & 0.0$\pm$0.0                             & 122.7$\pm$11.0                  & 0.0$\pm$0.0                             & 136.4$\pm$5.2                   & 0.0$\pm$0.0                             & 133.4$\pm$5.7                   & 0.0$\pm$0.0                             & 124.8$\pm$5.9                   \\
40                     & 10                     & 0.0$\pm$0.0                             & 22.0$\pm$35.1                   & 0.0$\pm$0.0                             & 52.4$\pm$103.5                  & 0.0$\pm$0.0                             & 66.7$\pm$45.7                   & 0.0$\pm$0.0                             & 64.9$\pm$50.2                   & 0.0$\pm$0.0                             & 33.5$\pm$24.3                   \\
240                    & 60                     & 0.0$\pm$0.0                             & 343.7$\pm$263.6                 & 0.0$\pm$0.0                             & 353.5$\pm$267.4                 & 0.0$\pm$0.0                             & 389.3$\pm$98.3                  & 0.0$\pm$0.0                             & 363.2$\pm$103.4                 & 0.0$\pm$0.0                             & 356.2$\pm$102.6                 \\
440                    & 110                    & 0.0$\pm$0.0                             & 1045.3$\pm$1484.5               & 0.0$\pm$0.0                             & 1052.5$\pm$1538.6               & 0.4$\pm$0.5                             & 1295.1$\pm$654.8                & 0.5$\pm$0.6                             & 1395.7$\pm$863.1                & 0.0$\pm$0.0                             & 1055.5$\pm$585.1                \\
640                    & 160                    & 0.0$\pm$0.0                             & 929.8$\pm$471.6                 & 0.0$\pm$0.0                             & 970.1$\pm$497.2                 & 0.0$\pm$0.0                             & 1040.7$\pm$169.1                & 0.0$\pm$0.0                             & 1055.6$\pm$200.5                & 0.0$\pm$0.0                             & 997.8$\pm$197.8                 \\
840                    & 210                    & 2.3$\pm$5.1                             & 3600.7$\pm$2048.6               & 3.5$\pm$8.0                             & 3520.7$\pm$2017.3               & 1.8$\pm$1.5                             & 3504.8$\pm$649.4                & 1.9$\pm$1.8                             & 3470.2$\pm$779.2                & 2.1$\pm$1.9                             & 3424.3$\pm$777.4                \\
1040                   & 260                    & 4.0$\pm$6.6                             & 5939.9$\pm$2770.5               & 3.3$\pm$7.4                             & 5767.3$\pm$2731.1               & 6.2$\pm$2.8                             & 5844.4$\pm$869.7                & 7.4$\pm$3.9                             & 5838.7$\pm$1033.4               & 6.3$\pm$3.8                             & 5654.2$\pm$1009.4               \\ \hline
\end{tabular}
\end{adjustbox}
}
{The values are averages over instance-behavior configurations with the same variable size, with a 95\% confidence interval.}
\end{table}

\begin{table}[]
\TABLE
{Runtime (seconds) of the {\MCMC} Hit-and-Run in the second step of the {\SAA} method. \label{tab:details-res-gen-eval}}
{
\centering
\begin{adjustbox}{width=0.85\textwidth}
\centering
\begin{tabular}{cc|c|c|c|c|c}
\hline
$n_x$ & $n_y$ & $\pi_{\text{Ntr}}(\cdot)$ & $\pi_{\text{Prp}}(\cdot)$ & $\pi_{\text{Frg}}(\cdot)$ & $\pi_{\text{Frg}^+}(\cdot)$ & \multicolumn{1}{c}{$\pi_{\text{Str}}(\cdot)$} \\ \hline
55   & 10  & 16.6$\pm$1.2     & 16.2$\pm$1.5     & 15.9$\pm$0.5     & 15.9$\pm$0.6     & 16.3$\pm$0.6     \\
55   & 15  & 19.3$\pm$2.4     & 19.8$\pm$1.5     & 20.1$\pm$0.5     & 19.9$\pm$0.6     & 19.2$\pm$0.7     \\
110  & 10  & 16.8$\pm$1.8     & 18.0$\pm$nan     & 17.5$\pm$nan     & 17.3$\pm$nan     & 17.8$\pm$nan     \\
165  & 5   & 11.6$\pm$0.9     & 11.9$\pm$0.9     & 11.8$\pm$0.4     & 12.0$\pm$0.4     & 12.0$\pm$0.4     \\
180  & 10  & 29.7$\pm$1.3     & 30.1$\pm$1.2     & 30.4$\pm$0.4     & 30.9$\pm$0.6     & 31.0$\pm$0.6     \\
225  & 5   & 18.4$\pm$0.9     & 18.4$\pm$0.9     & 18.9$\pm$0.3     & 18.6$\pm$0.3     & 18.9$\pm$0.4     \\
40   & 10  & 7.3$\pm$2.0      & 7.3$\pm$2.0      & 7.4$\pm$0.6      & 7.3$\pm$0.7      & 7.4$\pm$0.7      \\
240  & 60  & 93.5$\pm$46.2    & 92.8$\pm$45.9    & 92.7$\pm$14.7    & 93.6$\pm$17.9    & 86.6$\pm$16.5    \\
440  & 110 & 434.1$\pm$23.2   & 435.5$\pm$22.4   & 444.5$\pm$7.6    & 444.5$\pm$10.7   & 432.5$\pm$10.1   \\
640  & 160 & 804.6$\pm$206.3  & 802.9$\pm$205.7  & 813.7$\pm$66.8   & 825.5$\pm$80.7   & 828.3$\pm$81.5   \\
840  & 210 & 1774.1$\pm$226.3 & 1656.2$\pm$86.8  & 1474.7$\pm$112.2 & 1462.5$\pm$126.8 & 1626.7$\pm$29.2  \\
1040 & 260 & 2220.4$\pm$840.4 & 2175.7$\pm$826.9 & 2341.7$\pm$344.4 & 2354.2$\pm$420.1 & 2291.8$\pm$338.8 \\\hline
\end{tabular}
\end{adjustbox}
}
{The values are averages over instance-behavior configurations with the same variable size, with a 95\% confidence interval.}
\end{table}

\begin{figure}
    \FIGURE
    {
    \subcaptionbox{Average time (seconds). \label{subfig:comp-de-sw-time}}
    {\includegraphics[width=0.46\textwidth]{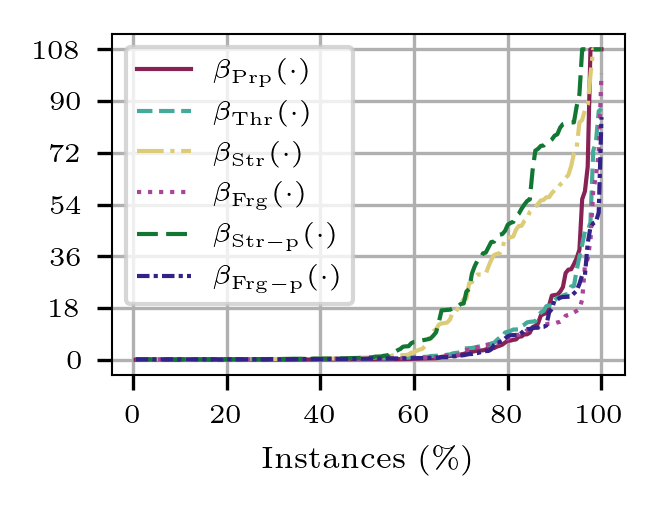}}
    \hfill
    \subcaptionbox{Average optimality gap [$\text{gap}_{opt}]$($\%$). \label{subfig:comp-de-sw-gap}}
    {\includegraphics[width=0.46\textwidth]{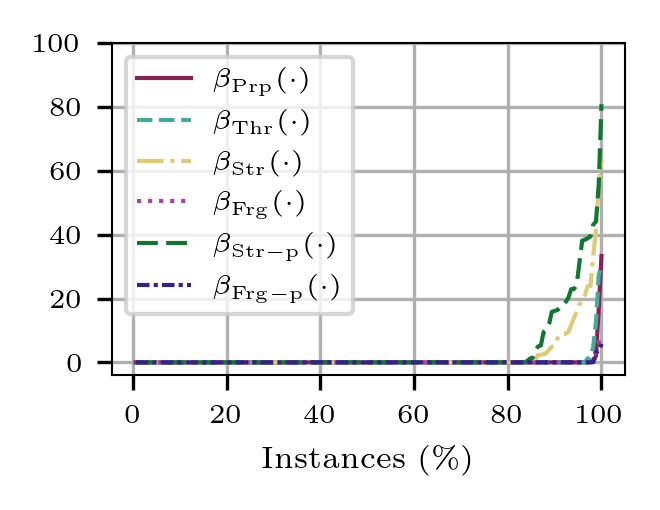}}
    }
    {Computational performance of the {\DE} program by strong-weak fixed follower behavior type. \label{fig:comp-de-sw}}
    {The values are averaged across parameterized variants of the same follower behavior (configurations), computed per instance.}
\end{figure}

\subsection{Comparison of Follower's Behavior}~\label{app:additional-results-comp}

Similarly to Section~\ref{subsubsec:objective-values}, in this section we also compare the follower optimal objective value $\varphi(x)$ corresponding to the optimal or best leader solutions obtained by the methods described in Section~\ref{subsubsec:performance-all-measures}, for different follower behaviors (i.e., endogenous measures), within the three-hour time limit. We use the optimal follower objective value under the optimistic approach as a reference. For a measure $\mu$, the follower optimality gap is defined as $\text{gap}_\varphi^\mu = 100 \% \cdot \boldsymbol{(}\varphi(\overline{x}^{\beta^1}) - \varphi(x^\mu)\boldsymbol{)}/ \lvert \varphi(\overline{x}_{\beta^1}) \rvert$, where $\overline{x}^{\beta^1}$ is the optimal leader solution obtained by the {\DE} program under the optimistic approach, and $x^\mu$ denotes the leader solution associated with measure $\mu$, namely $\overline{x}^\mu$ for the {\DE} program with strong-weak fixed measures, and $\overline{x}_N^\mu$ for the {\SAA} program with strong-weak and generalized decision-dependent measures, where the number of scenarios $N$ is defined in Section~\ref{subsec:comparison-behaviors}. Figure~\ref{fig:opt-follower-reference-gap} presents the median and quantiles of the follower optimal reference objective gap for each behavior, considering all instances, which follows a similar pattern to the follower expected reference objective gap.

\begin{figure}
    \FIGURE
    {
    {\includegraphics[height=4cm]{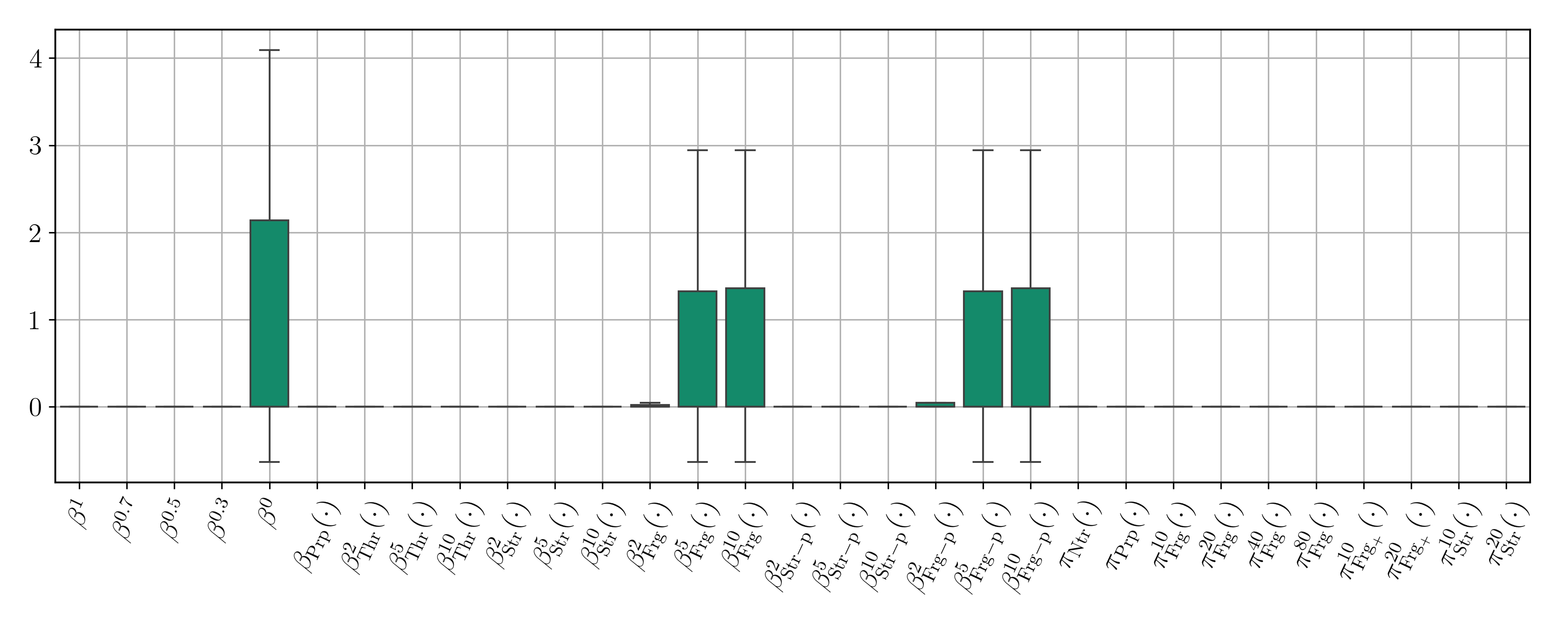}}
    }
    {Median and quantiles of the follower optimal reference objective gap [$\text{gap}_\varphi^\mu$] ($\%$) [All instances]. \label{fig:opt-follower-reference-gap}}
    {}
\end{figure}

\end{APPENDICES}

\end{document}